\documentclass{lmcs}
\pdfoutput=1

\usepackage{lastpage}
\lmcsdoi{20}{2}{19}
\lmcsheading{}{\pageref{LastPage}}{}{}%
{May~26,~2023}{Jun.~26,~2024}{}


\usepackage{amssymb}
\usepackage[latin1]{inputenc}
\usepackage{amsmath}
\usepackage{tikz-cd}
\usepackage{amsthm}
\usepackage[latin1]{inputenc}
\usepackage{graphicx}
\usepackage{hyperref}

\newtheoremstyle{maintheoremstyle}%
  {6pt}
  {6pt}
  {\itshape}
  {}
  {\bfseries}
  {{\bfseries .}}
  {5pt plus 1pt minus 1pt}
  {\thmname{#1} \thmnumber{#2}\thmnote{\normalfont{ (#3)}}}

\theoremstyle{maintheoremstyle}
\newtheorem{maintheorem}{Theorem}

\newenvironment{mainthm}[1]{\setcounter{maintheorem}{#1}\addtocounter{maintheorem}{-1}\begin{maintheorem}}{\end{maintheorem}}

\begin{document}
	
	\title[Robust non-computability of dynamical systems]{Robust non-computability of dynamical systems
		and computability of robust dynamical systems}
	\thanks{\textbf{Acknowledgments.} D. Gra\c{c}a was partially funded by
		FCT/MCTES through national funds and when applicable co-funded by EU
		funds under the project UIDB/50008/2020.
		\includegraphics[width=4.5mm]{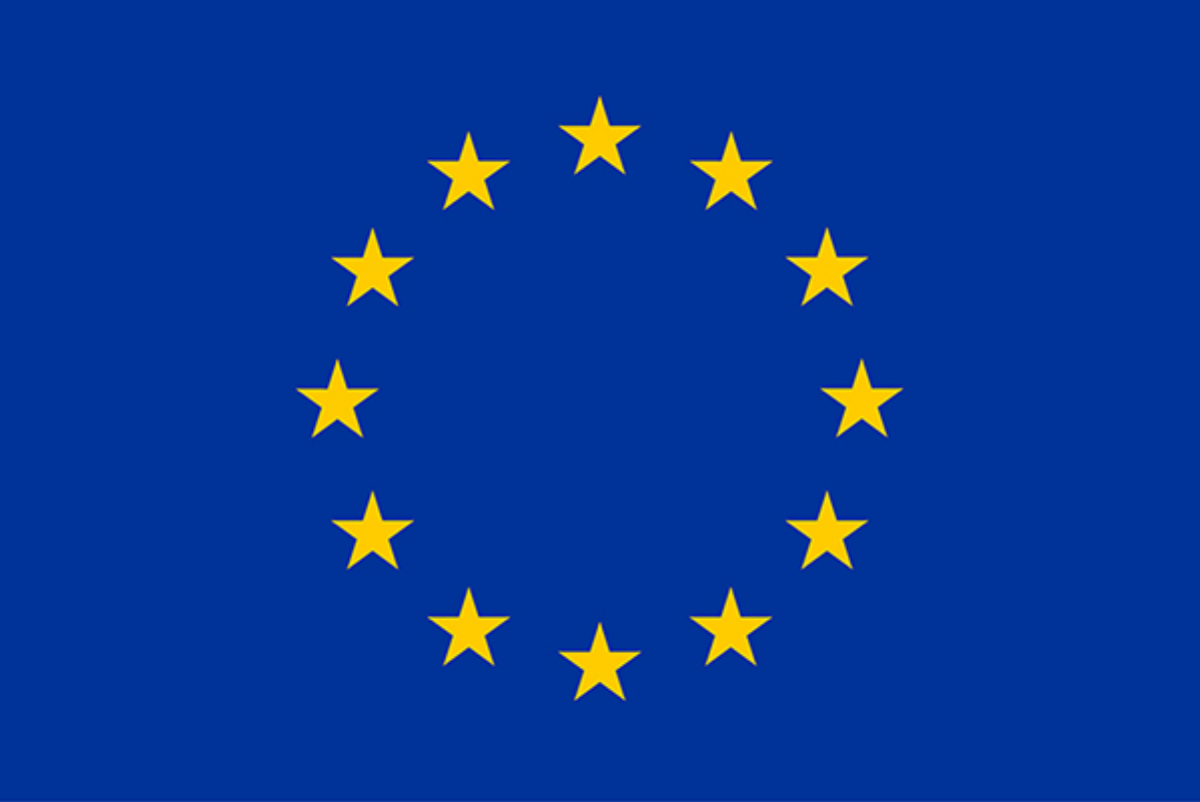} This project has received
		funding from the European Union's Horizon 2020 research and
		innovation programme under the Marie Sk\l {}odowska-Curie grant
		agreement No 731143.}
	
	\author[D. S. Gra\c{c}a]{Daniel S. Gra\c{c}a\lmcsorcid{0000-0002-0330-833X}}[a,b]
	\author[N. Zhong]{Ning Zhong}[c]
	
	\address{Universidade do Algarve, C. Gambelas, 8005-139 Faro, Portugal}
	
	\address{Instituto de Telecomunica\c{c}\~{o}es, 1049-001 Lisbon, Portugal}	
	
	\address{DMS, University of Cincinnati, Cincinnati, OH
		45221-0025, U.S.A.}

\begin{abstract}
In this paper, we examine the relationship between the stability of
the dynamical system $x^{\prime}=f(x)$ and the computability of its
basins of attraction. We present a computable $C^{\infty}$ system
$x^{\prime}=f(x)$ that possesses a computable and stable equilibrium
point, yet whose basin of attraction is robustly non-computable in a
neighborhood of $f$ in the sense that both the equilibrium point and
the non-computability of its associated basin of attraction persist
when $f$ is slightly perturbed. This indicates that local stability
near a stable equilibrium point alone is insufficient to guarantee
the computability of its basin of attraction. However, we also
demonstrate that the basins of attraction associated with a
structurally stable - globally stable (robust) - planar system defined on a compact set are
computable. Our findings suggest that the global stability of a
system and the compactness of the domain play a pivotal role in determining the computability of its
basins of attraction.
\end{abstract}

\maketitle

\section{Introduction} \label{introduction}

The focus of this paper is on examining the relationship between the
stability of the dynamical system
\begin{equation} \label{eq_main}
\frac{dx}{dt}=f(x)
\end{equation}
and the feasibility of computing the basin of attraction of a
(hyperbolic) equilibrium point.

The problem of computing the basin of attraction of an equilibrium
point can be viewed as a continuous variation of the discrete
Halting problem. In this paper, we will demonstrate that basins of
attraction can exhibit \emph{robust non-computability} for
computable systems. Specifically, we will present a computable
system represented by Equation \eqref{eq_main}  and a neighborhood
surrounding function $f$ which have the following properties: (i)
Equation \eqref{eq_main} has a computable equilibrium point, say
$s_f$, and the basin of attraction of $s_f$ is non-computable; (ii)
there are infinitely many computable functions within this
neighborhood; and (iii) for each and every computable function $g$
in this neighborhood,  the system described by $x^{\prime}=g(x)$
possesses a computable equilibrium point (near $s_f$) whose basin of
attraction is also non-computable. To the best of our knowledge,
this is the first instance where a continuous problem is
demonstrated to possess robust non-computability.

Equilibrium solutions, also known as equilibrium points or critical
points, correspond to the zeros of $f$ in \eqref{eq_main} and play a
vital role in dynamical systems theory. They are points where the
system comes to rest and are useful in determining the stability of
the system. By analyzing the system's behavior in the vicinity of an
equilibrium point, we can ascertain whether nearby trajectories
(i.e.~solutions of \eqref{eq_main}) will remain near that point
(stable) or move away from it (unstable).

\begin{figure}
    \begin{center}
    \includegraphics[width=7cm]{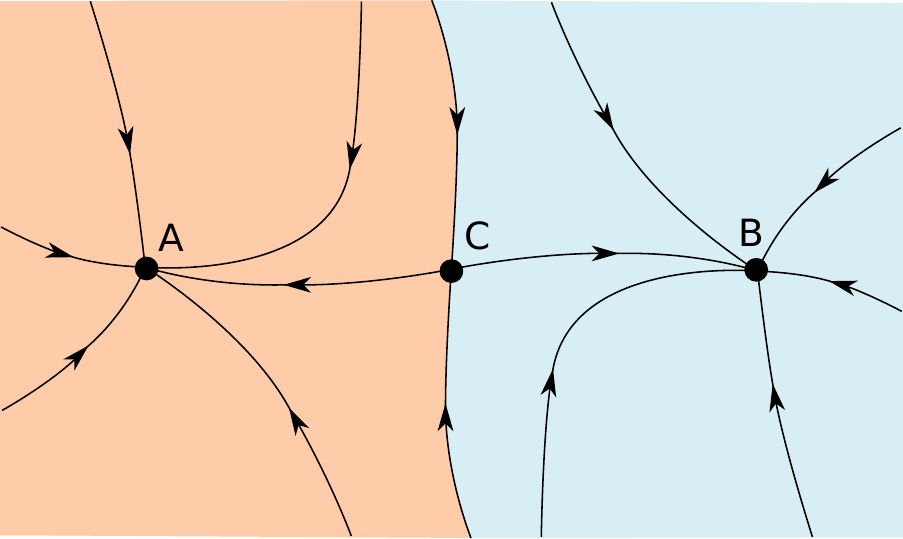}
    \caption{Example of a dynamical system having three equilibrium points $A,B,C$.
 The points $A$ and $B$ are sinks (i.e.~stable equilibrium points)
while $C$ is not (it is a so-called saddle equilibrium point). The
region in orange is the basin of attraction of $A$ while the region
in blue is the basin of attraction of $B$.}
    \label{Fig:basin}
    \end{center}
\end{figure}

The basins of attraction, on the other hand, represent the
collection of initial conditions with the property that their
associated trajectories converge to the corresponding equilibrium
point. This is pictured in Figure \ref{Fig:basin}. Thus, by
identifying the basins of attraction, we can predict the system's
long-term behavior for different initial conditions. This
information is essential in understanding and characterizing the
system's behavior, particularly in the context of complex systems. We also note that a basin of attraction is an open subset of some Euclidean space.

A sink of (\ref{eq_main}) is a special type of equilibrium point
where the system in the neighborhood of the equilibrium point is
well-behaved and stable. Here ``stable" refers to at least two
properties. First, each sink $s$ has a neighborhood $U$ with the
property that any trajectory that enters $U$ stays there and
converges exponentially fast to $s$ (this means that
$\|\phi_t(x)-s\|\leq \varepsilon e^{-\alpha t}$ for some
$\varepsilon,\alpha>0$, where $\phi_t(x)$ denotes the solution of
(\ref{eq_main}) at time $t\geq 0$ with initial condition
$\phi_0(x)=x\in U$. See  \cite[Theorem 1 on p.~130]{Per01}). Second,
the system is stable in the sense that if we replace $f$ in
(\ref{eq_main}) by a nearby function $\overline{f}$ then it will
continue to have a unique sink $\overline{s}$ in $U$ ($\overline{s}$
depends continuously on $\overline{f}$). In particular, when
$\overline{f}=f$ one has $s=\overline{s}$. See \cite[Theorem 1 on
p.~321]{Per01}. Moreover, trajectories of the new system will
behave (near the sink) similarly to the trajectories of the original
system (\ref{eq_main}) and will converge exponentially fast to
$\overline{s}$.

This means that if the system (\ref{eq_main}) is slightly perturbed
from a sink, it will eventually return to that point. In other
words, a sink point is robust locally under small perturbations.
Moreover, even if the dynamics of the system is (slightly)
perturbed, nearby trajectories will behave similarly to the original
system, providing a better understanding of the long-term behavior
of the system. This is particularly important in the study of
complex systems, where the stability of the system can be difficult
to determine analytically. The concept of robustness also allows for
the development of numerical methods for the study of dynamical
systems, which are crucial in many applications where analytical
methods are not feasible.

The widespread use of numerical algorithms in the analysis of
dynamical systems has made it crucial to determine which sets
associated with a system can be computed, and which ones cannot. In
essence, a set is computable if it can be accurately plotted or
numerically described to any desired degree of precision.
Equilibrium points and basins of attraction are examples of such
sets.

Several studies \cite{Zho09}, \cite{GZ15} revealed that the basin of
attraction of a sink may not be computable, even if the system is
analytic and computable, and the sink is computable. Furthermore,
non-computability results about dynamical systems are not restricted only to basins of
attraction for differential equations (some examples can be found in e.g.~\cite{PR79},
\cite{PR81}, \cite{PZ97}, \cite{BY06}, \cite{GZB07}, \cite{HS08a}, \cite{GBC09}, \cite{GHR11}, \cite{GHR12}, \cite{GZB12}, \cite{CRY18}, \cite{RY20}, \cite{BP20}, \cite{GHRS20}, \cite{CMP21}, \cite{CMPP21}, \cite{CFHR22}). These discoveries
highlight the need to understand the limitations of numerical
methods in the analysis of dynamical systems. In particular, it
raises the question of whether non-computability results are
``typical" or if they represent ``exceptional" scenarios that are
unlikely to have practical significance. In this paper, we
specifically concentrate on investigating the non-computability of
basins of attractions, as this phenomenon can be viewed as a
continuous-time counterpart to the halting problem.

Moreover, it's worth noting that numerical computations have finite
precision, and hyperbolic sinks are robust under small
perturbations. As a result, it's worth considering if the
non-computability remains under small perturbations. If it does not,
the non-computability may be ignored in physical realities. In this
paper, we show that the non-computability in computing the basin of
attraction cannot be overlooked in the sense that the
non-computability is robust under small perturbations. The following
is our first main result (the precise statement is presented in
section~3).

\begin{mainthm}{1}
There exists a computable $C^{\infty}$
function $f$ for which the system (\ref{eq_main}) possesses a
computable sink $s_0$, but the basin of attraction of $s_0$ is
non-computable. Moreover, this non-computability is robust and
persists under small perturbations.
\end{mainthm}

It is worth noting that Theorem A establishes that local stability
in the vicinity of a sink is insufficient to guarantee the
computability of the basin of attraction at the sink.

We also provide a discrete-time variant of this theorem. Actually, this result will be proved first and then used to prove Theorem A.

\begin{mainthm}{2}
There exists an analytic and computable
function $f$ for which the discrete-time dynamical system defined by the iteration of $f$ possesses a
computable sink $s_0$, but the basin of attraction of $s_0$ is
non-computable. Moreover, this non-computability is robust and
persists under small perturbations.
\end{mainthm}

The third main theorem of the paper provides an answer to the
question of which dynamical systems have computable basins of
attraction in the plane $\mathbb{R}^2$. The precise statement of the following theorem is presented in section~4. This theorem applies to planar structurally stable system. Intuitively, a system is structurally stable if perturbing ``a little bit'' the dynamics $f$ of \eqref{eq_main} will not change the qualitative shape of the dynamics (see Definition \ref{Def:structurally-stable} for a formal definition). For this reason we will sometimes say that a system is globally stable if it is structurally stable.

\begin{mainthm}{3}
The map that links each structurally stable
planar system defined on a compact set to the set of basins of attraction of its sinks is
computable.
\end{mainthm}

This theorem provides a positive result that complements the
non-computability result presented in Theorem A. It
implies that for the large class of structurally stable - globally
stable - planar systems, it is possible to numerically compute the
basins of attraction (as open sets) of their equilibrium points and periodic
orbits. It is worth noting that the set of structurally stable
planar systems defined on a compact disk $K$ forms an open and dense
subset of the set of all planar ($C^1$-) dynamical systems defined
on $K$.

Taken together, Theorems A and C demonstrate that global stability
is a crucial element in determining the feasibility of numerically
computing the basins of attraction of dynamical systems, at least in
the case of ordinary differential equations.

\section{Preliminaries}

\label{Sec:prelims}

\subsection{Computable analysis}

\label{Subsec:comptanalysis}

Let $\mathbb{N}, \mathbb{Z}, \mathbb{Q}$, and $\mathbb{R}$ be the
set of non-negative integers, integers, rational numbers, and real
numbers, respectively. Assuming familiarity with the concept of
computable functions defined on $\mathbb{N}$ with values in
$\mathbb{Q}$, we note that there exist several distinct but equally
valid approaches to computable analysis, dating back to the work of
Grzegorczyk and Lacombe in the 1950s. For the purposes of this
paper, we adopt the oracle Turing machine version presented in e.g.~\cite{Ko91}.

\begin{defi} \label{def_oracle} A rational-valued function
$\phi: \mathbb{N}\to \mathbb{Q}$ is called an oracle for a real
number $x$ if it satisfies $|\phi(m) - x| < 2^{-m}$ for all $m$.
\end{defi}

\begin{defi} \label{def_function} Let $S$ be a subset of $\mathbb{R}$, and let $f:
S\to \mathbb{R}$ be a real-valued function on $S$. Then $f$ is said
to be computable if there is an oracle Turing Machine $M^{\phi}(n)$
such that the following holds: If $\phi$ is an oracle for $x\in
S$, then for every $n\in \mathbb{N}$, $M^{\phi}(n)$ returns a
rational number $q$ such that $|q-f(x)|<2^{-n}$.
\end{defi}

The definition can be extended to functions defined on a subset of
$\mathbb{R}^d$ with values in $\mathbb{R}^l$.

\begin{defi} \label{def_open_set} Let $U$ be a bounded open subset of
$\mathbb{R}^d$. Then $U$ is called computable if there are
computable functions $a, b: \mathbb{N}\to \mathbb{Q}^d$ and $r, s:
\mathbb{N}\to \mathbb{Q}$ such that the following holds: $U=\cup
_{n=0}^{\infty}B(a(n), r(n))$ and $\{ \overline{B(b(n), s(n))}\}_n$
lists all closed rational balls in $\mathbb{R}^d$ which are disjoint
from $U$, where $B(a, r)=\{ x\in \mathbb{R}^d \, : \, |x-a| < r\}$
is the open ball in $\mathbb{R}^d$ centered at $a$ with the radius
$r$ and $\overline{B(a, r)}$ is the closure of $B(a, r)$.
\end{defi}

By definition, a planar computable bounded open set can be rendered
on a computer screen with arbitrary magnification. A closed subset
$K$ of $\mathbb{R}^d$ is considered computable if its complement
$\mathbb{R}^d\setminus K$ is a computable open subset of
$\mathbb{R}^d$, or equivalently, if the distance function $d_{K}:
\mathbb{R}^d \to \mathbb{R}$ defined as $d_{K}(x) = \inf _{y\in D}
\| y - x\|$ is computable.

The concept of Turing computability can be extended to encompass a
broader range of function spaces and the maps that operate on them.
The definitions \ref{def_function} and \ref{def_open_set} indicate
that an object is deemed (Turing) computable if it can be
approximated with arbitrary precision through computer-generated
approximations. Formalizing this idea to carry out computations on
infinite objects such as real numbers, we encode those objects as
infinite sequences of rational numbers (or equivalently, sequences
of any finite or countable set $\Sigma$ of symbols), using
representations (see \cite{Wei00} for a complete development). A
represented space is a pair $(X; \delta)$ where $X$ is a set,
$\delta$ is a coding system (or naming system) on $X$ with codes
from $\Sigma$ having the property that
$\mbox{dom}(\delta)\subseteq\Sigma^{\mathbb{N}}$ and $\delta:
\Sigma^{\mathbb{N}}\to X$ is an onto map. Every
$q\in\mbox{dom}(\delta)$ satisfying $\delta(q)=x$ is called a
$\delta$-name of $x$ (or a name of $x$ when $\delta$ is clear from
context). Naturally, an element $x\in X$ is computable if it has a
computable name in $\Sigma^{\mathbb{N}}$. The notion of
computability on $\Sigma ^{\mathbb{N}}$ is well established, and
$\delta$ lifts computations on $X$ to computations on
$\Sigma^{\mathbb{N}}$. The representation $\delta$ also induces a
topology $\tau_{\delta}$ on $X$, where $\tau_{\delta} = \{
U\subseteq X: \, \delta^{-1}(U)\text{ is open in }\operatorname{dom}
(\delta)\}$ is called the final topology of $\delta$ on $X$.

The notion of computable maps between represented spaces now arises
naturally. A map $\Phi: (X;\delta_{X})\to(Y;\delta_{Y})$ between two
represented spaces is computable if there is a computable map $\phi:
\Sigma^{\mathbb{N}}\to\Sigma^{\mathbb{N}}$ such that
$\Phi\circ\delta_{X}=\delta_{Y}\circ\phi$  as depicted below (see
e.g.~\cite{BHW08}).
\begin{center}
    \begin{tikzcd}

        \Sigma^{\mathbb{N}} \arrow[r, "\phi"] \arrow[d, "\delta_{X}"]
        & \Sigma^{\mathbb{N}} \arrow[d, "\delta_{Y}"] \\
        X \arrow[r, "\Phi"]
        & Y
    \end{tikzcd}
\end{center}
Informally speaking, this means that there is a computer program
$\phi$ that outputs a name of $\Phi(x)$ when given a name of $x$ as
input. Since $\phi$ is computable, it transforms every computable
element in $\Sigma^{\mathbb{N}}$ to a computable element in
$\Sigma^{\mathbb{N}}$.
Another fact about computable maps is that computable maps are
continuous with respect to the corresponding final topologies
induced by $\delta_X$ and $\delta_Y$.

\subsection{Dynamical systems}

\label{Subsec:DynamicalSystem}

Discrete-time dynamical systems are defined by the
iteration of a map $g:\mathbb{R}^{d}\rightarrow\mathbb{R}^{d}$,
while continuous-time systems are defined by an ordinary
differential equation (ODE) of the form $x^{\prime}=f(x)$, where $f:\mathbb{R}^{d}\rightarrow\mathbb{R}^{d}$.
Regardless of the type of system, the notion of trajectory is
fundamental. In the discrete-time case, a trajectory starting at the
point $x_{0}$ is defined by the sequence of iterates of $g$ as
follows
\[
x_{0},g(x_{0}),g(g(x_{0})),\ldots,g^{[k]}(x_0),\ldots
\]
where $g^{[k]}$ denotes the $k$th iterate of $g$, while in the
continuous time case it is the solution, a function $\phi (f,
x_0)(\cdot)$ of time $t$,  to the following
initial-value problem
\[
\left\{
\begin{array}
[c]{l}
x^{\prime}=f(x)\\
x(0)=x_{0}
\end{array}
\right.
\]

In the realm of dynamical systems, a set $A$ is considered forward
invariant if any trajectory starting on $A$ remains on $A$
indefinitely for any positive time. If an invariant set consists of
only one point, it is called an equilibrium point. For a dynamical
system defined by \eqref{eq_main}, an equilibrium point must be a
zero of $f$. Similarly, for a discrete-time dynamical system defined
by $g$, an equilibrium point must be a fixed point of $g$ (i.e. it
satisfies $g(x)=x$) or, equivalently, it must be a zero of $g(x)-x$.

If trajectories nearby an invariant set converge to this set, then
the invariant set is called an attractor. The basin of attraction
for a given attractor $A$ is the set of all points
$x\in\mathbb{R}^{d}$ such that the trajectory starting at $x$
converges to $A$ as $t\rightarrow\infty$. Attractors come in
different types, including points, periodic orbits, and strange
attractors. Equilibrium points are the simplest type of attractor.

An equilibrium point $x_0$ of \eqref{eq_main} is \emph{hyperbolic}
if none of the eigenvalues of the Jacobian matrix $Df(x_0)$ have
zero real part. In particular, if all the eigenvalues of $Df(x_0)$
have a negative real part, then we are have a \emph{sink}. A sink
has all the properties mentioned in Section~\ref{introduction}. In
particular given a sink $s$ there is a neighborhood $U$ such that
any trajectory starting in $U$ stays there and converges
exponentially fast to $s$. If an hyperbolic equilibrium point is not
a sink, then given any neighborhood of this point, there will be a
trajectory that will never reach this equilibrium point.

A similar approach can be applied to discrete-time dynamical
systems. Specifically, an equilibrium point $x_0$ of the
discrete-time dynamical system defined by $g$ is hyperbolic if  none of the
eigenvalues of $Dg(x_0)$ belong to the unit circle. On the other
hand, an equilibrium point $x_0$ is considered a sink if all the
eigenvalues of $Dg(x_0)$ have an absolute value less than 1.

We will now discuss the concept of ($C^1$-)perturbations. First, we
will introduce some notations. Let $C^k(A; \mathbb{R}^l)$ denote the
set of all $k$-times continuously differentiable functions from a
subset $A$ of $\mathbb{R}^d$ to $\mathbb{R}^l$. If $l=d$, we simply
write $C^k(A)$ for $C^k(A; \mathbb{R}^d)$. Suppose $W$ is an open
subset of $\mathbb{R}^d$ and $f: W \to \mathbb{R}^d$ is a $C^1$
vector field. In the field of dynamical systems and differential
equations, a perturbation of $f$ is another $C^1$ vector field $g: W
\to \mathbb{R}^d$ that is ``$C^1$-close to $f$". To be more precise:

\begin{defi} \label{def_norm} Let $f\in C(W)$ (resp. $f\in C^1(W)$), the $C$-norm of $f$ is
defined to be $\|f\|=\sup_{x\in W}\|f(x)\|$ (resp.~the $C^1$-norm of
$f$ is defined to be $\|f\|_1 = \sup_{x\in W}\|f(x)\| + \sup_{x\in
W}\| Df(x)\|$), where  $\|\cdot\|= $ denotes the max-norm on
$\mathbb{R}^d$ or the usual norm of the matrix
$Df(x)$, depending on the context.
\end{defi}

Note that for $x\in \mathbb{R}^d$, the max-norm is given by
$\|x\|=\max_{1\leq i\leq d}|x_i|$. It is possible that $\|
f\|_1=\infty$ if the number is unbounded. The $C^1$-norm $\|\cdot
\|_1$ has many of the same formal properties as norms for vector
spaces. For $\epsilon>0$, an $\epsilon$-neighborhood of $f$ in
$C^1(W)$ is defined as the set $\{ g\in C^1(W): \| g - f\|_1 <
\epsilon\}$. Any function $g$ in this neighborhood is called an
$\epsilon$-perturbation of $f$.

\begin{rem}\label{Rem:many-functions}
    Upon observation, it can be noted that for any function $f: W
\rightarrow \mathbb{R}^l$, if $f$ is computable with a finite
$\|f\|_1$, then in any $\epsilon$-neighborhood (in $C^1$-norm) $\mathcal{N}$, there
exist infinitely many computable $C^1$ functions which are distinct
from $f$. For example, $f_{\alpha}, \bar{f}_{\alpha}, \tilde{f}_{\alpha}\in \mathcal{N}$ for any rational $\alpha$ satisfying $0<\alpha<\epsilon$, where (the operations are done componentwise) $f_{\alpha}(x)=f(x)+\alpha$, $\bar{f}_{\alpha}(x)=f(x)+\alpha\sin x$, $\tilde{f}_{\alpha}(x)=f(x)+e^{-\alpha(1+\|x\|^2)}$.
\end{rem}

\begin{figure}
	\begin{center}
		\includegraphics[width=12cm]{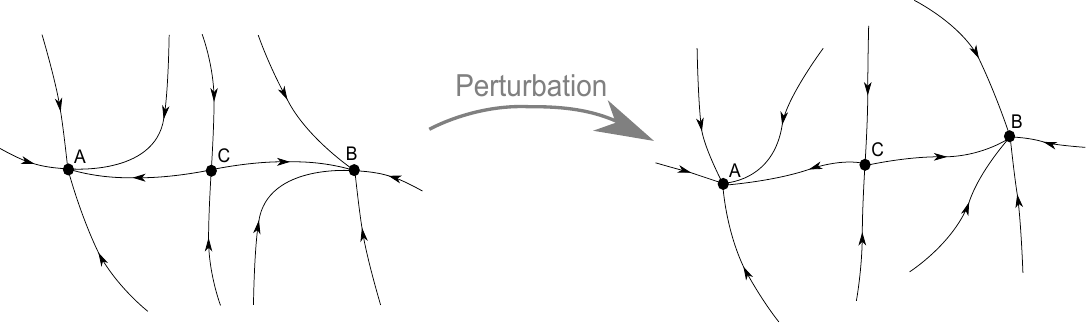}
	\end{center}
	\caption{An example of a structurally stable system on the left. Even if perturbed the main properties of the system persist. For example, there is a connection between the sink $A$ and the saddle $C$ and similarly for $B$ and $C$ which persists under (small perturbation).}\label{fig:structurally-stable}
\end{figure}

\begin{figure}
	\begin{center}
		\includegraphics[width=105mm]{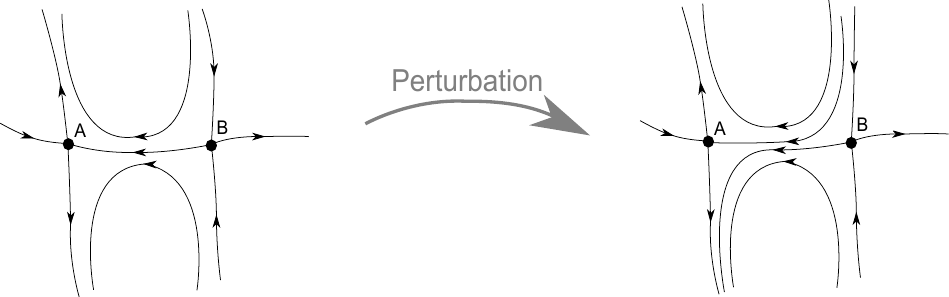}
	\end{center}
	\caption{An example of a structurally unstable system on the left, known as a saddle connection. We can find a perturbation, which can be assumed to be as small as we want, that is able to break the connection between the saddles $A$ and $B$.}\label{fig:structurally-unstable}
\end{figure}

Next we present the notion of structural stability (see Figures \ref{fig:structurally-stable} and \ref{fig:structurally-unstable} for a picture).
\begin{defi}\label{Def:structurally-stable}
	A planar dynamical system $dx/dt = f(x)$,
	where $f\in \mathcal{V}(K)$, is structurally stable if there
	exists some $\varepsilon>0$ such that for all $g\in C^{1}(K)$
	satisfying $\left\Vert f-g\right\Vert _{1}\leq\varepsilon$, the
	trajectories of $dy/dt = g(y)$ are homeomorphic to the trajectories
	of $dx/dt = f(x)$. In other words, there exists a homeomorphism $h:K\to K$
	such that if $\gamma$ is a oriented trajectory of $dx/dt = f(x)$, then
	$h(\gamma)$ is a oriented trajectory of $dy/dt = g(y)$.
\end{defi}

\section{Proof of Theorem B -- robust non-computability in the discrete-time case}

\label{sec:EverywhereNoncomputable}

In this section, we provide an example demonstrating the existence
of a computable and analytic function $f: \mathbb{R}^3 \rightarrow
\mathbb{R}^3$ that defines a discrete-time dynamical system
satisfying the following conditions:
\begin{enumerate}[(i)]
\item $f$ has a hyperbolic sink $s$ that is computable.
\item The basin of attraction of $s$ is non-computable.
\item There exists a neighborhood $\mathcal{N}$ (in $C^1$-norm) of $f$
such that for every function $g \in \mathcal{N}$, $g$ has a
hyperbolic sink $s_g$ that is computable from $g$, and the basin of
attraction of $s_g$ is non-computable.
\end{enumerate}

The construction demonstrates that the non-computability of
computing the basins of attraction can remain robust under small
perturbations and sustained throughout an entire neighborhood.

It is worth noting that the function $f$ inherits strong
computability from its analyticity, which implies that every order
derivative of $f$ is computable. Furthermore, in any
$C^1$-neighborhood of $f$, there exist infinitely many computable
functions (see Remark \ref{Rem:many-functions}).

We will make use of the following example that is explicitly constructed in \cite[Section 4]{GZ15}.

\begin{propC}[\cite{GZ15}]\label{Pro:sink-theorem}
	There is an analytic and computable function $f: \mathbb{R}^3 \to
	\mathbb{R}^3$ with the following properties:
\begin{enumerate}[(a)]
\item the restriction, $f_{M}: \mathbb{N}^3 \to \mathbb{N}^3$, of $f$ on
$\mathbb{N}^3$ is the transition function of a one-tape universal Turing
machine $M$, where each configuration of $M$ is coded as an element
of $\mathbb{N}^3$ (see below for an exact description of the
coding). Without loss of generality, $M$ can be assumed to have just
one halting configuration; e.g. just before ending, clear the tape and go to a unique halting state;
thus the halting configuration $s\in\mathbb{N}^3$ is unique. We also assume that $f_M$ is defined over $s$ and that $f(s)=s$.
\item the halting configuration $s$ of $M$ is a computable
sink of the discrete-time evolution of $f$.
\item the basin of attraction of $s$ is non-computable.
\item there exists a constant $\lambda \in (0, 1)$ such that if $x_0$
is a configuration of $M$, then for any $x\in \mathbb{R}^3$,
\begin{equation} \label{f-contraction}
\| x - x_0\| \leq 1/4 \quad \Longrightarrow \quad \| f(x) -
f(x_0)\|\leq \lambda \| x - x_0\|
\end{equation}
\end{enumerate}
\end{propC}

The coding  used in Proposition \ref{Pro:sink-theorem} of the one-tape configuration of the Turing machine $M$ with tape contents $\ldots B B B a_{-k}\ldots a_{-1} a_0 a_1 \ldots a_n B B B \ldots$ and state $q\in\{1,\ldots,m\}$, where $B$ is the blank symbol, $a_0$ is the symbol being read by the tape head, and $a_i$ are symbols of an alphabet $\Gamma$ with, without loss of generality, at most 10 symbols including the blank symbol, is given by $x=(x_1,x_2,x_3)\in\mathbb{N}^3$ where
\begin{align*}
x_{1}  & =\iota(a_{0})+\iota(a_{1})10+\ldots+\iota(a_{n})10^{n}\\
x_{2}  & =\iota(a_{-1})+\iota(a_{-2})10+\ldots+\iota(a_{-k})10^{k-1}\\
x_{3}  & =q
\end{align*}
and $\iota:\Gamma\to\{0,1,\ldots,9\}$ is an injective function (coding) such that $\iota(B)=0$.

Roughly, the proof of this result relies on the use of interpolation techniques to extend the transition function from $\mathbb{N}^3$ to $\mathbb{R}^3$. In particular one can use trigonometric interpolation to obtain a function $\omega:\mathbb{R}\to\mathbb{R}$ such that $\omega(i)=i$ for $i\in\{0,1,\ldots,9\}$ which ``recovers'' from $x_1$ the symbol $a_0$ being read by the head since $\omega(x_1)=\iota(a_0)$. Then noting that the transition function of $M$ is only defined over a finite number of pairs $(q,a_0)$, we can extend this transition function to $\mathbb{R}$ using polynomial interpolation and then update the value of the triple $x=(x_1,x_2,x_3)\in\mathbb{N}^3$ coding the configuration. This idea serves as the foundation for establishing property (a) as discussed in the preceding study [GCB08, p. 333], wherein property (a) is substantiated alongside an exploration of additional robustness characteristics.

Nevertheless, in \cite[p. 333]{GCB08}, the Halting configuration is linked not to a sink but  to a set of values within a ball. In the work presented in \cite{GZ15}, this construction is enhanced to ensure the Halting configuration is explicitly associated with a sink $s$, thereby affirming property (b). This accomplishment is realized by guaranteeing the satisfaction of property (d). Given the undecidability of the Halting problem for a universal Turing machine, it consequently establishes property (c).

Indeed, condition (d) can be strengthened in two ways: (i) the contraction rate $\lambda$ in \eqref{f-contraction} can be selected to be arbitrarily small, and (ii) robustness against perturbations in the dynamics can be attained, mirroring the approach employed in \cite[Theorem 1]{GCB08}.

\begin{thm}\label{Th:contraction}
	For every $\lambda\in (0,1)$, $\varepsilon\in (0,1/4]$, and every transition function $f_{M}: \mathbb{N}^3 \to \mathbb{N}^3$ of a one-tape Turing machine $M$, there is an extension $f:\mathbb{R}^3 \to \mathbb{R}^3$ of $f_M$ with the following properties, where $x_{0}\in\mathbb{N}^{3}$ denotes a configuration of $M$:
	
	\begin{enumerate}
		\item It satisfies condition \eqref{f-contraction}.
		
		\item If $\left\Vert f-g\right\Vert \leq(1-\lambda)\varepsilon<1/4$, then for all
		$j\in\mathbb{N}$ one has
		\[
		\left\Vert x-x_{0}\right\Vert \leq\varepsilon\text{ \ \ }\Rightarrow\text{
			\ \ }\left\Vert g^{[j]}(x)-f_{M}^{[j]}(x_{0})\right\Vert \leq\varepsilon.
		\]
		
		\item If $M$ has a unique halting configuration $s$ and $f_{M}(s)=s$, similarly to condition (a) of Proposition \ref{Pro:sink-theorem}, then $s$ is a computable sink of the discrete-time evolution of $f$.
		
	\end{enumerate}
\end{thm}

Note that this theorem implies that, by considering a universal Turing machine with an assumed unique halting configuration $s \in \mathbb{N}^3$ as mentioned earlier, 
we can conclude that for any selected $\lambda\in (0,1)$, there exists a function $f:\mathbb{R}^3 \to \mathbb{R}^3$ satisfying properties (a), (b), (c) and property (d) for the specific value of $\lambda$.

\begin{proof}[Proof of Theorem \ref{Th:contraction}]
	Let's commence by demonstrating condition 1. We first note that the function $\sigma:\mathbb{R\rightarrow R}$ defined as $\sigma(x)=x-0.2\sin(2\pi x)$ is a uniform contraction around integers (see \cite[Proposition 5]{GCB08}), i.e.~it satisfies the following property
	\[
	\left\vert x-n\right\vert \leq1/4\text{ \ \ }\Longrightarrow\text{
		\ \ }\left\vert \sigma(x)-n\right\vert <\lambda_{1/4}\left\vert
	x-n\right\vert
	\]
	where $\lambda_{1/4}=0.4\pi-1\approx0.256637$ for any $n\in\mathbb{Z}$.
	Moreover, given a one-tape Turing machine $M$ one can find some $\bar{\lambda
	}\in(0,1)$ such that the transition function $f_{M}:\mathbb{N}^{3}
	\rightarrow\mathbb{N}^{3}$ admits an extension $\bar{f}:\mathbb{R}
	^{3}\rightarrow\mathbb{R}^{3}$ with the property that (see \cite[Theorem 4]{GZ15})
	\[
	\| x - x_0\| \leq 1/4 \quad \Longrightarrow \quad \| \bar{f}(x) -
	\bar{f}(x_0)\|\leq \bar{\lambda} \| x - x_0\|
	\]
	(recall that $x_0\in\mathbb{N}^{3}$ is the coding of a configuration). Now given some arbitrary $\lambda
	\in(0,1)$ one can find some $k\in\mathbb{N}$ such that $0<\lambda_{1/4}
	^{k}\bar{\lambda}<\lambda<1$. Hence, by applying $\sigma^{\lbrack k]}$ to each
	component of $\bar{f}$, we get a function $f:\mathbb{R}^{3}\rightarrow
	\mathbb{R}^{3}$ which satisfies \eqref{f-contraction}. This also implies condition 3, i.e.~that $s$ is a sink.
	
	Concerning condition 2, we note that for $j=1$ one has
	\begin{align*}
		\left\Vert g(x)-f_{M}(x_{0})\right\Vert  & =\left\Vert g(x)-f(x_{0}
		)\right\Vert \\
		& \leq\left\Vert g(x)-f(x)\right\Vert +\left\Vert f(x)-f(x_{0})\right\Vert \\
		& \leq(1-\lambda)\varepsilon+\lambda\left\Vert x-x_{0}\right\Vert \\
		& \leq \varepsilon 	.
	\end{align*}
	Proceeding by induction for $j>1$ we conclude the desired result.
\end{proof}

In the remaining of this section, the symbols $f$ and $s$ are
reserved for this particular function and its particular sink for a universal Turing machine $M$ whose
transition function is $f_M$.

We now show in the following that there is a $C^1$-neighborhood
$\mathcal{N}$ of $f$ -- computable from $f$ and $Df(s)$ -- such that
for every $g\in \mathcal{N}$, $g$ has a sink $s_{g}$ -- computable
from $g$ -- and the basin of attraction of $s_{g}$ is
non-computable.

The following proposition is a computable version of a classical result of dynamical systems theory. See, for example, the Proposition in \cite[p. 305]{HS74}. While the original version of this result in \cite{HS74} does not refer to computability, it is straightforward to check that its proof also proves the computable version of this proposition. The proof of this
proposition has nothing to do with differential equations; rather, it
depends on the invertibility of $Dh(z_0)$.

\begin{prop}\label{Prop:pertubed-sink}
	Suppose that $h:\mathbb{R}^{n}\rightarrow\mathbb{R}^{n}$ is of class $C^{1}$
	and suppose that $z_{0}\in\mathbb{R}^{n}$ is such that $Dh(z_{0})$ is
	invertible. Then one can compute from $h$ and $z_{0}$ rationals $\epsilon
	,\delta>0$ such that for the neighborhood $U=B(z_{0},\epsilon)$ of $z_{0}$ and
	for any $C^{1}$ function $g:\mathbb{R}^{n}\rightarrow\mathbb{R}^{n}$
	satisfying $\left\Vert g-h \right\Vert _{1}<\delta$ one has:
	
	\begin{enumerate}
		\item $g$ is injective in $U$;
		
		\item $h(z_{0})\in g(U)$;
		
		\item $Dg$ is invertible in $U$;
		
		\item For any rational $\bar{\epsilon}>0$, with $\bar{\epsilon}\leq\epsilon$,
		one can choose a rational $\bar{\delta}>0$, with $\delta\geq\bar{\delta}$, such that if $\left\Vert
		g-h \right\Vert _{1}<\mathcal{\bar{\delta}}$, then there is a $\bar{z}_{0}\in
		B(z_{0},\bar{\epsilon})$ satisfying $g(\bar{z}_{0})=h(z_{0})$.
	\end{enumerate}
\end{prop}

\begin{cor}
	Suppose that $s$ is a sink of a function $f:\mathbb{R}^{n}\rightarrow
	\mathbb{R}^{n}$. Then one can compute from $f$, $Df$, and $s$ a neighborhood $U$ of
	$s$ and a $C^{1}$-neighborhood $\mathcal{N}$ of $f$ such that for any
	$g\in\mathcal{N}$, $g$ has a unique sink $s_{g}$ in $U$ which is computable from $f$, $Df$, $s$, $g$, and $Dg$. Moreover, for any rational
	$\epsilon>0$ one can choose $\mathcal{N}$ so that $\Vert s_{g}-s\Vert
	<\epsilon$.
\end{cor}
\begin{proof}
	Immediate from the previous proposition by taking $h=f$ and $z_0=s$.
\end{proof}

\begin{mainthm}{2}
There is a $C^{1}$-neighborhood $\mathcal{N}$ of $f$ (computable from $f$ and $Df(s)$) such that for any $g\in\mathcal{N}$, $g$ has a sink $s_{g}$ (computable from $g$) and the basin of attraction $W_{g}$ of $s_{g}$ is non-computable.
\end{mainthm}

\begin{proof}
	We first note that if $s$ is a sink of $f$, then $Df(s)-I$ is invertible,
	i.e.~$\det(Df(s)-I)\neq0$. Indeed, if $\det(Df(s)-I)=0$ then there would be a
	non-zero vector $v$ such that $(Df(s)-I)v=0$, a contradiction to the
	assumption that $s$ is a sink and thus that all eigenvalues of $Df(s)$ are
	less than 1 in absolute value. This implies that the sink $s$ is a zero of
	$\bar{f}(x)=f(x)-x$ and that $D\bar{f}(s)$ is invertible.
	
Specifically, if we consider the function $f$ as described in Proposition \ref{Pro:sink-theorem} on the set $K=\{x\in\mathbb{R}^{n}:\left\Vert x-s\right\Vert \leq1/4\}$, it becomes evident that it possesses a unique fixed point within this set. This solitary fixed point is a sink, a consequence of the contraction property given in \eqref{f-contraction}.
Furthermore, due to Theorem \ref{Th:contraction}, we can assume that $\lambda=1/4$ and that the condition 2 of
	that theorem holds with $\varepsilon=1/8$. This implies that if $g\in
	C^{1}(K)$ is such that $\left\Vert f-g\right\Vert _{1}\leq1/8$, then for any
	configuration $x_{0}\in\mathbb{N}^{3}$ of $M$, and any $x\in\overline
	{B(x_{0},1/4)}$, we have the following estimate:
	\begin{align}
		&  \Vert g(x)-g(x_{0})\Vert\nonumber\\
		&  \leq\Vert(g-f)(x)-(g-f)(x_{0})\Vert+\Vert f(x)-f(x_{0})\Vert\nonumber\\
		&  \leq\Vert D(g-f)\Vert\,\Vert x-x_{0}\Vert+\lambda\Vert x-x_{0}
		\Vert\nonumber\\
		&  \leq(1/8+\lambda)\,\Vert x-x_{0}\Vert\nonumber\\
		&  \leq \frac{3}{8}\Vert x-x_{0}\Vert\label{g-contraction}
	\end{align}

	Since $3/8<1$, it follows that $g$ is a contraction in $\overline{B(x_{0},1/4)}$ for
	every configuration $x_{0}$ of $M$. 
 
Now, let's utilize Proposition \ref{Prop:pertubed-sink} with the function $h=\bar{f}$ and  $z_{0}=s$. We will compute positive rational values $\epsilon$ and $\delta$ that fulfill the conditions stated in the proposition. We also take $\bar{\epsilon}
	=\min(\epsilon/2,1/16)$ for condition 4 of that proposition and obtain the
	respective $\bar{\delta}>0$ which we can assume without loss of generality to
	satisfy $\bar{\delta}\leq1/8$. Given some $\bar{g}\in C^{1}(K)$ satisfying
	$\left\Vert \bar{g}-\bar{f}\right\Vert _{1}\leq\bar{\delta}$, let $s_{\bar{g}
	}$ denote the unique zero of $\bar{g}$ in $B(s,\epsilon)$, which corresponds
	to a sink $s_{g}$ of $g(x)=\bar{g}(x)+x$
	satisfying $\left\Vert s_{g}-s\right\Vert \leq1/16$. Note that, in this case,
	\[
	\left\Vert g(x)-f(x)\right\Vert =\left\Vert g(x)-x-(f(x)-x)\right\Vert
	=\left\Vert \bar{g}(x)-\bar{f}(x)\right\Vert .
	\]
	We now show that for any $g\in C^{1}(K)$ satisfying $\left\Vert g-f\right\Vert
	_{1}\leq\bar{\delta}$ and for any configuration $x_{0}\in\mathbb{N}^{3}$ of
	$M$, $M$ halts on $x_{0}$ if and only if $x_{0}\in W_{g}$, where $W_{g}$
	denotes the basin of attraction of $s_{g}$.
	
	First we assume that $x_{0}\in W_{g}$. Then, by definition of basin of
	attraction of a sink, $g^{[j]}(x_{0})\rightarrow s_{g}$ as $j\rightarrow
	\infty$. Hence, there exists $n\in\mathbb{N}$ such that $\Vert g^{[n]}
	(x_{0})-s_{g}\Vert<\frac{1}{16}$, which in turn implies that
	\begin{align*}
		&  \Vert f_{M}^{[n]}(x_{0})-s\Vert\\
		&  \leq\Vert f_{M}^{[n]}(x_{0})-g^{[n]}(x_{0})\Vert+\Vert g^{[n]}(x_{0}
		)-s_{g}\Vert+\Vert s_{g}-s\Vert\\
		&  \leq\frac{1}{8}+\frac{1}{16}+\frac{1}{16}\leq\frac{1}{8}.
	\end{align*}
	Due to \eqref{f-contraction} and because $s$ is a sink of $f$, it follows that $f_{M}
	^{[n]}(x_{0})=s$. Hence, $M$ halts on $x_{0}$ and, moreover,
there
	exists $n\in\mathbb{N}$ such that $f_{M}^{[j]}(x_{0})=s$ for all $j\geq n$.
	Then for all $j\geq n$, it follows from Theorem \ref{Th:contraction}
	that
	\begin{align*}
		&  \Vert g^{[j]}(x_{0})-s\Vert\\
		&  \leq\Vert g^{[j]}(x_{0})-f_{M}^{[j]}(x_{0})\Vert+\Vert f_{M}^{[j]}
		(x_{0})-s\Vert\\
		&  =\Vert g^{[j]}(x_{0})-f_{M}^{[j]}(x_{0})\Vert\leq1/8
	\end{align*}
	The inequality implies that $\{g^{[j]}(x_{0})\}_{j\geq n}\subset
	\overline{B(s,1/8)}$. Because $s_{g}$ is a sink of $g$ satisfying $\Vert
	s-s_{g}\Vert<\frac{1}{16}$, it follows that $g(s_{g})=s_{g}$ and $s_{g}
	\in\overline{B(s,\bar{\epsilon})}\subset\overline{B(s,1/4)}$. Since $s$ is a
	configuration of $M$ -- the halting configuration of $M$ -- it follows from
	(\ref{g-contraction}) that $g$ is a contraction on $\overline{B(s,1/4)}$.
	Thus, $\Vert g^{[n+j]}(x_{0})-s_{g}\Vert=\Vert g^{[n+j]}(x_{0})-g^{[n+j]}
	(s_{g})\Vert\leq(\theta_{\lambda})^{j}\Vert g^{[n]}(x_{0})-s_{g}
	\Vert\rightarrow0$ as $j\rightarrow\infty$. Consequently, $g^{[j]}
	(x_{0})\rightarrow s_{g}$ as $j\rightarrow\infty$, This implies that $x_{0}\in
	W_{g}$.
	
	To prove that $W_{g}$ is non-computable, the following stronger inclusion is
	needed: if $M$ halts on $x_{0}\in\mathbb{N}^3$, then $B(x_{0},1/8)\subset W_{g}$. Consider
	any $x\in\overline{B(x_{0},1/8)}$. Since $x_{0}\in W_{g}$ and $g$ is a
	contraction on $\overline{B(x_{0},1/8)}$ due to (\ref{g-contraction}), it
	follows that
	\[
	\Vert g^{[j]}(x)-g^{[j]}(x_{0})\Vert\leq(3/8)^{j}\Vert x-x_{0}\Vert
	\rightarrow0\quad\mbox{as $j\to \infty$}
	\]
	Since $x_{0}\in W_{g}$, $g^{[j]}(x_{0})\rightarrow s_{g}$ as $j\rightarrow
	\infty$. Hence, $g^{[j]}(x)\rightarrow s_{g}$ as $j\rightarrow\infty$. This
	implies that $x\in W_{g}$. Moreover, if $M$ does not halt on $x_{0}$, then
	$\overline{B(x_{0},1/8)}\cap W_{g}=\varnothing$ due to Theorem \ref{Th:contraction}.
	
	It remains to show that $W_{g}$ is non-computable. Suppose otherwise that
	$W_{g}$ was computable. We first note that $W_{g}=\bigcup_{t\in\mathbb{N}}
	\phi_{-t}(B(x_{0},\epsilon))$ is an open set (we recall that $\phi_t(x)$ denotes the solution of
	(\ref{eq_main}) at time $t\in\mathbb{R}$ with initial condition
	$\phi_0(x)=x\in U$) since $\phi_{t}$ is continuous
	for every $t\in\mathbb{R}$ (this is a well-known fact that follows from the
	formula \eqref{error-initial-condition}) and furthermore $\phi_{t}^{-1}
	=\phi_{-t}$. Then the distance function $d_{\mathbb{R}^{3}\setminus W_{g}}$ is
	computable. We can use this computability to solve the halting problem.
	Consider any initial configuration $x_{0}\in\mathbb{N}^{3}$, and compute
	$d_{\mathbb{R}^{3}\setminus W_{g}}(x_{0})$. If $d_{\mathbb{R}^{3}\setminus
		W_{g}}(x_{0})>\frac{1}{9}$ or $d_{\mathbb{R}^{3}\setminus W_{g}}(x_{0}
	)<\frac{1}{8}$, halt the computation. Since $\epsilon>0$, this computation
	always halts.
	
	Now we use the fact that either $\overline{B(x_{0},1/8)}$ is fully contained
	in $W_{g}$ or otherwise $\overline{B(x_{0},1/8)}$ does not intersect $W_{g}$
	and is thus fully contained in $\mathbb{R}^{3}\setminus W_{g}$. If
	$d_{\mathbb{R}^{3}\setminus W_{g}}(x_{0})>\frac{1}{9}>0$, then $\overline
	{B(x_{0},1/8)}$ is not fully contained in $\mathbb{R}^{3}\setminus W_{g}$
	which implies that $x_{0}\in W_{g}$, or equivalently, the Turing machine $M$
	halts on $x_{0}$ Otherwise, if $d_{\mathbb{R}^{3}\setminus W_{g}}(x_{0}
	)<\frac{1}{8}$, then there are points of $\overline{B(x_{0},1/8)}$ in
	$\mathbb{R}^{3}\setminus W_{g}$ and this can only happen if $\overline
	{B(x_{0},1/8)}\subseteq\mathbb{R}^{3}\setminus W_{g}$, which implies that
	$x_{0}\not \in W_{g}$, or equivalently, $M$ does not halt on $x_{0}$.
	Therefore, if $W_{g}$ was computable, then we could solve the halting problem,
	which is a contradiction. Hence, we conclude that $W_{g}$ is non-computable.
\end{proof}

\begin{rem}Theorem B demonstrates that
non-computability can maintain its strength when considering
standard topological structures, as in the study of natural
phenomena such as identifying invariant sets of a dynamical system.
This robustness can manifest in a powerful way: the
non-computability of the basins of attraction persists continuously
for every function that is ``$C^1$ close to $f$''.
\end{rem}

\section{Proof of Theorem A -- robust non-computability\texorpdfstring{\\}{} in the continuous-time case}

In the previous section, we demonstrated that a discrete-time
dynamical system defined by the iteration of a map, say
$\bar{f}:\mathbb{R}^{3}\rightarrow\mathbb{R}^{3}$, has a computable
sink with a non-computable basin of attraction, and that this
non-computability property is robust to perturbations. In this
section, we extend this result to continuous-time dynamical systems.
Specifically, we prove the existence of a computable $C^{\infty}$
map $f:\mathbb{R}^{7}\rightarrow \mathbb{R}^{7}$ such that the ODE
$y^{\prime}=f(y)$ has a computable sink with a non-computable basin
of attraction. Moreover, this non-computability property is robust
to small perturbations in $f$.

To be more precise, we show that there exists some $\varepsilon>0$
such that if $g:\mathbb{R}^{7}\rightarrow \mathbb{R}^{7}$ is another
$C^{\infty}$ map with $\left\Vert f-g\right\Vert
_{1}\leq\varepsilon$, then the ODE $y^{\prime}=g(y)$ also has a sink
(computable from $g$ and located near the sink of $y^{\prime}=f(y)$)
with a non-computable basin of attraction. This means that the
non-computability of the basin of attraction is a robust property of
the underlying dynamical system.

Overall, this result shows that the non-computability of basin of
attraction is not limited to discrete-time dynamical systems, but is
also present in continuous-time dynamical systems, and is a robust
property that persists under small perturbations.

To obtain this result, we will employ a technique that involves
iterating the map $\bar{f}$ with an ODE. This
technique has been explored in several previous papers, including
\cite{Bra95}, \cite{CMC00}, \cite{CMC02}, \cite{GCB08}, and \cite{GZ23}. However, we need additional requirements which are not ensured by the original technique, namely we need to ensure that the
resulting ODE still has a computable sink and that the
non-computability property is robust to perturbations.
Hence, it is imperative to expound upon the details of the previous constructions, as our primary approach revolves around the progressive enhancement of preceding constructions, aiming to acquire additional properties essential for the substantiation of Theorem A.

\subsection{Iterating a map with an ODE}

The basic idea
to iterate a map with an ODE is to start with a \textquotedblleft
targeting\textquotedblright\ equation with
the format
\begin{equation}
    x^{\prime}=c(b-x)^{3}\phi(t)\label{Eq:targeting_ODE}
\end{equation}
where $b$ is the \emph{target value} and
$\phi:\mathbb{R}\mathbb{\rightarrow R}$ is a continuous function
which satisfies $\int_{t_{0}}^{t_{1}}\phi(t)dt>0$ and $\phi(t)\geq0$
over $[t_{0},t_{1}]$. This is a separable ODE which can be
explicitly solved. Using the solution one can show that for any
$\gamma>0$ (the value $\gamma$ is called the \emph{targeting error}
for reasons which will be clear in a moment), if one chooses
\begin{equation}
    c\geq\frac{1}{2\gamma^{2}\int_{t_{0}}^{t_{1}}\phi(t)dt}\label{Eq:c}
\end{equation}
in (\ref{Eq:targeting_ODE}), then $\left\vert x(t)-b\right\vert
<\gamma$ for all $t\geq t_{1}$, independent of the initial condition
$x(t_{0})$. Note also that if $\phi(t)=0$ for all $t\in\lbrack
t_{0},t_{1}]$, then $x(t)=x(t_{0})$ for all $t\in\lbrack
t_{0},t_{1}]$. This targeting equation is
the basic construction block for iterating a map $\tilde{f}
:\mathbb{R\rightarrow R}$, which extends a corresponding function
$\tilde {f}_{\mathbb{N}}:\mathbb{N\rightarrow N}$.

To iterate $\tilde{f}$ (with an ODE) we pick $t_{1}-t_{0}=1/2$, a
continuous periodic function $\phi:\mathbb{R\rightarrow R}$ of
period 1, which satisfies $\phi(t)\geq0$ for
$t\in]0,1/2[$, $\phi(t)=0$ for $t\in \lbrack1/2,1]$, and
$\int_{0}^{1}\phi(t)dt>0$, a constant $c$ satisfying (\ref{Eq:c})
with $\gamma=1/4$, and a $C^{\infty}$ function
$r:\mathbb{R\rightarrow R}$ with the property that
$r(k+\varepsilon)=k$ for all $k\in\mathbb{Z}$ and all
$0\leq\varepsilon\leq1/4$ (i.e.~$r$ returns the integer part of its
argument $x$ whenever $x$ is within distance $\leq1/4$ of an
integer). Although the exact expressions of $\phi$ and $r$
are irrelevant to the construction, it is worth noticing that choices can be made (see e.g.~\cite[p.~344]
{GCB08}, replacing $\theta_{j}$ in (20) of that paper by the function $\chi$
given by (\ref{Eq:smooth-Heaviside}) below) so that $\phi$ and $r$ are
$C^{\infty}$.

Then the ODE
\begin{equation}
    \left\{
    \begin{array}
        [c]{l}
        z_{1}^{\prime}=c(\tilde{f}(r(z_{2}))-z_{1})^{3}\phi(t)\\
        z_{2}^{\prime}=c(r(z_{1})-z_{2})^{3}\phi(t+1/2)
    \end{array}
    \right.  \label{Eq:iteration}
\end{equation}
will iterate $\tilde{f}$ in the sense that the continuous flow
generated by (\ref{Eq:iteration}) starting near any integer value
will stay close to the (discrete) orbit of $\tilde{f}$, as we will
now see.
Suppose that at the initial time
$t=0$, we have $\left\vert z_{1}(0)-x_{0}\right\vert \leq1/4$ and
$\left\vert z_{2}(0)-x_{0}\right\vert \leq1/4$ for some
$x_{0}\in\mathbb{Z}$. During the first half-unit interval $[0,1/2]$,
we have $\phi(t+1/2)=0$, and thus $z_{2}^{\prime}(t)=0$.
Consequently, $z_{2}(t)=z_{2}(0)$, and hence $r(z_{2})=x_{0}$.
Therefore, the first equation of (\ref{Eq:iteration}) becomes a
targeting equation (\ref{Eq:targeting_ODE}) on the interval
$[0,1/2]$ where the target is
$\tilde{f}(r(z_{2}))=\tilde{f}(x_{0})$. Thus, we have $\left\vert
z_{1}(1/2)-\tilde{f}(x_{0})\right\vert \leq1/4$.

In the next half-unit interval $[1/2,1]$, the behavior of $z_{1}$
and $z_{2}$ switches. We have $\phi(t)=0$, and thus
$z_{1}(t)=z_{1}(1/2)$, which implies that
$r(z_{1})=\tilde{f}(x_{0})$. Hence, the second equation of
(\ref{Eq:iteration}) becomes a targeting equation
(\ref{Eq:targeting_ODE}) on the interval $[0,1/2]$ where the target
is $r(z_{1})=\tilde{f}(x_{0})$. Thus, we have $\left\vert
z_{2}(1)-\tilde{f}(x_{0})\right\vert \leq1/4$.

In the next unit interval $[1,2]$, the same behavior repeats itself,
and therefore we conclude that we have $\left\vert
z_{1}(2)-\tilde{f}(\tilde{f}(x_{0}))\right\vert \leq1/4$ and
$\left\vert z_{2}(2)-\tilde{f}(\tilde{f}(x_{0}))\right\vert
\leq1/4$. In general, for any $k\in\mathbb{N}$ and $t\in\lbrack
k,k+1/2]$, we will have $\left\vert
z_{1}(k)-\tilde{f}^{[k]}(x_{0})\right\vert \leq1/4$, $\left\vert
z_{2}(k)-\tilde{f}^{[k]}(x_{0})\right\vert \leq1/4$, and $\left\vert
z_{2}(t)-\tilde{f}^{[k]}(x_{0})\right\vert \leq1/4$. In other words,
the flow of (\ref{Eq:iteration}) starting near any integer value
stays close to the orbit of $\tilde{f}$.

Notice also that by choosing $\gamma=1/8$ instead of $\gamma=1/4$,
we can make (\ref{Eq:iteration}) robust to perturbations of
magnitude $\leq1/8$, since under these conditions the system
\begin{equation}
    \left\{
    \begin{array}
        [c]{l}
        \bar{z}_{1}^{\prime}=c(\tilde{f}(r(\bar{z}_{2}))-\bar{z}_{1})^{3}\phi
        (t)+\xi_{1}(t)\\
        \bar{z}_{2}^{\prime}=c(r(\bar{z}_{1})-\bar{z}_{2})^{3}\phi(t+1/2)+\xi_{2}(t)
    \end{array}
    \right.  \label{Eq:robust-iteration}
\end{equation}
still satisfies $\left\vert z_{1}(k)-\tilde{f}
^{[k]}(x_{0})\right\vert \leq1/4$, $\left\vert z_{2}(k)-\tilde{f}^{[k]}
(x_{0})\right\vert \leq1/4$, and $\left\vert z_{2}(t)-\tilde{f}^{[k]}
(x_{0})\right\vert \leq1/4$ for all $k\in\mathbb{N}$ and
$t\in\lbrack k,k+1/2]$, where $\left\vert \xi_{1}(t)\right\vert
\leq1/8$, $\left\vert \xi _{2}(t)\right\vert \leq1/8$ for all
$t\in\mathbb{R}$, and $\left\vert z_{1}(0)-x_{0}\right\vert
\leq1/8$, $\left\vert z_{2}(0)-x_{0}\right\vert
\leq1/8$.  Indeed, in $[0,1/2]$ we have $\phi(t+1/2)=0$ and hence $\bar{z}
_{2}^{\prime}=\xi_{2}(t)$ which yields $\left\vert z_{2}(t)-z_{2}
(0)\right\vert \leq\int_{0}^{1/2}\left\vert \xi_{2}(t)\right\vert
dt\leq(1/2)(1/8)=1/16$ and thus $\left\vert z_{2}(t)-x_{0}\right\vert
\leq\left\vert z_{2}(t)-z_{2}(0)\right\vert +\left\vert z_{2}(0)-x_{0}
\right\vert \leq1/16+1/8=3/16$ for all $t\in\lbrack0,1/2]$. Therefore
$\tilde{f}(r(\bar{z}_{2}))=\tilde{f}(x_{0})$ in $[0,1/2]$. Using an analysis
similar to that performed in \cite[p.~346]{GCB08}, where the \textquotedblleft
perturbed\textquotedblright\ targeting ODE
\begin{equation}
    x^{\prime}=c(b-x)^{3}\phi(t)+\xi(t) \quad (\mbox{with $\left\vert \xi(t)\right\vert \leq\rho$})
    \label{Eq:target-perturbed}
\end{equation}
is studied, we conclude that if $c$ satisfies (\ref{Eq:c}), then
$\left\vert x(t_{1})-b\right\vert <\gamma+\rho\cdot (t_{1}-t_{0})$.
In the present case $t_{1}-t_{0}=1/2$ and $\rho=1/8$, and thus
$\left\vert z_{1}(1/2)-\tilde{f}(x_{0})\right\vert
\leq1/8+(1/8)(1/2)=3/16$. Similarly, since $\phi(t)=0$, on $[1/2,1]$
we conclude that $\left\vert z_{1}(t)-\tilde{f}(x_{0})\right\vert
\leq\left\vert z_{1}(1/2)-\tilde{f}(x_{0})\right\vert
+\int_{1/2}^{1}\left\vert \xi _{2}(t)\right\vert
dt\leq3/16+(1/2)(1/8)=1/4$ for all $t\in\lbrack1/2,1]$. Therefore
$r(\bar{z}_{1})=\tilde{f}(x_{0})$ in $[1/2,1]$ and thus $\left\vert
z_{2}(1)-\tilde{f}(x_{0})\right\vert \leq1/8+(1/8)(1/2)=3/16$. By
repeating this procedure on subsequent intervals, we conclude that
$\left\vert z_{1}(k)-\tilde{f}^{[k]}(x_{0})\right\vert \leq1/4$,
$\left\vert z_{2}(k)-\tilde{f}^{[k]}(x_{0})\right\vert \leq1/4$ and
$\left\vert z_{2}(t)-\tilde{f}^{[k]}(x_{0})\right\vert \leq1/4$ for
all $k\in\mathbb{N}$ and $t\in\lbrack k,k+1/2]$.

The above procedure can be readily extended to iterate (with an ODE)
the three-dimensional map
$\bar{f}:\mathbb{R}^{3}\rightarrow\mathbb{R}^{3}$ of the previous
section by assuming that
$\bar{f}=(\bar{f}_{1},\bar{f}_{2},\bar{f}_{3})$, where
$\bar{f}_{i}:\mathbb{R}^{3}\rightarrow\mathbb{R}$ is a component of
$\bar{f}$
for $i=1,2,3$. To accomplish this, it suffices to consider the ODE
\begin{equation}
    \left\{
    \begin{array}
        [c]{l}
        u_{1}^{\prime}=c(\bar{f}_{1}(r(v_{1}),r(v_{2}),r(v_{3}))-u_{1})^{3}\phi(t)\\
        u_{2}^{\prime}=c(\bar{f}_{2}(r(v_{1}),r(v_{2}),r(v_{3}))-u_{2})^{3}\phi(t)\\
        u_{3}^{\prime}=c(\bar{f}_{3}(r(v_{1}),r(v_{2}),r(v_{3}))-u_{3})^{3}\phi(t)\\
        v_{1}^{\prime}=c(r(u_{1})-v_{1})^{3}\phi(t+1/2)\\
        v_{2}^{\prime}=c(r(u_{2})-v_{2})^{3}\phi(t+1/2)\\
        v_{3}^{\prime}=c(r(u_{3})-v_{3})^{3}\phi(t+1/2)
    \end{array}
    \right.  \label{Eq:iteration-3D}
\end{equation}
This ODE works like (\ref{Eq:iteration}), but applies componentwise
to each component $\bar{f}_{1},\bar{f}_{2},\bar{f}_{3}$.

\subsection{Ensuring that the halting configuration is a sink}\label{Sec:sink}

We have so far presented the basic technique used in \cite{Bra95}, \cite{CMC00}, \cite{CMC02}, \cite{GCB08} (several improvements exist from paper to paper). However, this is not sufficient for the purposes of the present paper and a
few problems still need to be addressed in order to achieve our
desired results. Specifically, we must: 
\begin{enumerate}[(i)]
	\item  acquire an autonomous
	system of the form $y^{\prime}=f(y)$ rather than a non-autonomous
	one like (\ref{Eq:iteration-3D});
	
	\item demonstrate the existence of
	a sink with a non-computable basin of attraction;
	
	\item establish that both the sink and the non-computability of the basin
	of attraction are resilient to perturbations.
\end{enumerate}
 In this sense we need to improve the constructions from previous papers. 

In this subsection we will improve the construction of the previous subsection to address (i) and (ii), much along the lines of what is done in \cite{GZ15}, although it will be important to present all the details for when addressing (iii). The condition (iii) will be addressed in the next subsection.

To address problem (i), one possible solution would be to introduce
a new variable $z$ that satisfies $z^{\prime}=1$ and $z(0)=0$,
effectively replacing $t$ in (\ref{Eq:iteration-3D}) with $z$.
However, this approach would not be compatible with problem (ii)
because the component $z$ would grow infinitely and never converge
to a value, which is necessary for the existence of a sink.

One potential solution to this problem is to introduce a new
variable $z$ such that $z(0)=0$ and $z^{\prime}=1$ until the Turing
machine $M$ halts, and then set $z^{\prime}=-z$ afterwards so that
the dynamics of $z$ converge to the sink at $0$ in one-dimensional
dynamics. Since $z$ will replace $t$ as the argument of $\phi$ in
(\ref{Eq:iteration-3D}), we also need to modify $\phi$ such that
when $M$ halts, the components of $u=(u_{1},u_{2},u_{3})$ and
$v=(v_{1},v_{2},v_{3})$ still converge to a sink that corresponds to
the unique halting configuration of $M$.

In order to describe the dynamics of $z$,
we first need to introduce several auxiliary tools.
Consider the $C^{\infty}$ function $\chi$ defined by
\begin{equation}
    \chi(x)=\left\{
    \begin{array}
        [c]{ll}
        0 & \text{if }x\leq0\text{ or }x\geq1\\
        e^{\frac{1}{x(x-1)}}\text{ \ \ } & \text{if }0<x<1.
    \end{array}
    \right.  \label{Eq:smooth-Heaviside}
\end{equation}
Notice that $\chi$,  as well as all its derivatives, is computable.
Now consider the $C^{\infty}$ function $\zeta$ defined by
$\zeta(0)=0$ and
\[
\zeta^{\prime}(x)=c\chi(x)
\]
where $c=\left(  \int_{0}^{1}e^{\frac{1}{x(x-1)}}dx\right)  ^{-1}$, which is a
$C^{\infty}$ version of Heaviside's function (see also \cite[p.~4]{Cam02a})
since $\zeta(x)=0$ when $x\leq0$, $\zeta(x)=1$ when $x\geq1$, and
$0<\zeta(x)<1$ when $0<x<1$. Notice that $\zeta$ is computable since the
solution of an ODE with $C^1$ computable data is computable \cite{GZB07}. Similar properties are trivially obtained for the
function $\zeta_{a,b}$, where $a<b$, defined by
\[
\zeta_{a,b}(x)=\zeta\left(  \frac{x-a}{b-a}\right)  =\left\{
\begin{array}
    [c]{ll}
    0 & \text{if }x\leq a\\
    \ast & \text{if }a<x<b\\
    1\text{ \ \ } & \text{if }x\geq b
\end{array}
\right.
\]
where $\ast$\ is a value in $]0,1[$ that depends on $x$. Let us now
update the function $\phi$ to be used in (\ref{Eq:iteration-3D}).
Recall that, in the previous section, we introduced the map
$\bar{f}:\mathbb{R}^{3}\rightarrow\mathbb{R}^{3}$ ($\bar{f}$ is
called $f$ in the previous section), which simulates a Turing
machine by encoding each configuration as (an approximation of) a
triplet $(w_{1},w_{2},q)\in\mathbb{N}^{3}$ (for more details, see
\cite{GCB08}). Here, $w_{1}$ encodes the part of the tape to the
left of the tape head (excluding the infinite sequence of
consecutive blank symbols), $w_{2}$ encodes the part of the tape
from the location of the tape head up to the right, and $q$ encodes
the state. We typically assume that ${1,\ldots,m}$ encode the
states, and $m$ represents the halting state. In
(\ref{Eq:iteration-3D}), $v_{3}$ gives the current state of the
Turing machine $M$, i.e., $v_{3}(t)=q_{k}$ for all $t\in\lbrack
k,k+1/2]$ if the state of $M$ after $k$ steps is $q_{k}$.
Additionally, $v_{3}(t)\in\lbrack q_{k},q_{k+1}]$
($v_{3}(t)\in\lbrack q_{k+1},q_{k}]$) if $q_{k}\leq q_{k+1}$
($q_{k+1} < q_{k}$, respectively) and $t\in\lbrack k+1/2,k+1].$
Define
\begin{equation}
    \bar{\phi}(t,v_{3})=\phi(t)+\zeta_{m-1/4,m-3/16}(v_{3}).\label{Eq:phi_bar}
\end{equation}
We note that $\phi(x),\zeta_{a,b}(x)\in[0,1]$ for any
$x\in\mathbb{R}$. Moreover, if $M$ halts in $k$ steps, then
$\bar{\phi}(t,v_{3}(t))=\phi(t)$ for $t\leq k-1/2$, and
$1\leq\bar{\phi}(t,v_{3}(t))\leq2$ when $t\geq k$. Let us now
analyze what happens when $t\in[k-1/2,k]$. We observe that
$v_{3}(t)$ will increase in this interval from the value of
approximately $q_{k-1}$ until it reaches a $1/4$-vicinity of
$q_{k}=m$. Until that happens, $\bar{\phi}(t)=\phi(t)$. Once
$v_{3}(t)$ is in $[m-1/4,m-3/16]$, we get that
$\bar{\phi}(t,v_{3}(t))=\phi(t)+\zeta_{m-1/4,m-3/16}(v_{3}(t))>\phi(t)$,
and if we use $\bar{\phi}(t,v_{3}(t))$ instead of $\phi(t)$ in the
first three equations of (\ref{Eq:iteration-3D}), the respective
targeting equations still have the same dynamics but with a faster
speed of convergence. Thus, because the targeting error is
$\gamma=1/8$, at a certain time $t^{\ast}$ we will have
$v_{3}(t)\geq m-3/16$ for all $t\geq t^{\ast}$. From this point on,
we will have (note that $1\geq\phi(t)$)
\begin{equation}
    \bar{\phi}(t, v_3)=\phi(t)+1\geq1\geq\phi(t)\label{Eq:phi_1-halt}
\end{equation}
and thus all 6 equations of (\ref{Eq:iteration-3D}) will become
\textquotedblleft locked\textquotedblright\ with respect to their
convergence, regardless of the value of $\phi(t)$ (and
$\phi(t+1/2)$). In other words, for $t\geq t^{\ast}$, the
convergence of the 6 equations of (\ref{Eq:iteration-3D}) is
guaranteed even if $\phi(t)=0$ or $\phi(t+1/2)=0$ for all $t\geq
t^{\ast}$. This means that from this moment $t$ can take any value.
In particular, from that moment we can replace $t$ by a variable $z$
which converges to $0$, as desired from our considerations described
above.

Let
\begin{equation}
    z^{\prime}=1-\zeta_{m-3/16,m-1/8}(v_{3}(t))(z+1),\text{ \ \ }z(0)=0
    \label{Eq:new_time}
\end{equation}
Notice that $z^{\prime}=1$ for all $t\leq t^{\ast}$. Hence $z(t)=t$
for all $t\leq t^{\ast}$. Once $v_{3}(t)$ reaches the value $m-1/8$
at time $t^{\ast\ast}>t^{\ast}$, we have $v_{3}(t)\geq m-1/8$. Hence
$z^{\prime}=-z$ for all $t\geq t^{\ast\ast}$ and thus $z$ will
converge exponentially fast to 0 for $t\geq t^{\ast\ast}$.

Let us now show that $x_{halt}=(0,0,m,0,0,m,0)\in \mathbb{R}^{7}$ is
a sink (recall that $(w_{1},w_{2},q)\in\mathbb{N}^{3}$ encodes a
configuration when simulating the Turing machine with the map
$\bar{f}:\mathbb{R}^{3}\rightarrow\mathbb{R}^{3}$). We may assume
that the machine cleans its tape before halting, thus generating the
halting configuration $(0,0,m)\in\mathbb{N}^{3}$. First we should
note that, as pointed out in \cite[Section 5.5]{GZ15}, all 6
equations of
(\ref{Eq:iteration-3D}) are variations of the ODE
\[
z^{\prime}=-z^{3}
\]
which has an equilibrium point at $z=0$, but is not hyperbolic, and
thus $z=0$ cannot be a sink. Therefore $x_{halt}$ cannot be a sink
of (\ref{Eq:iteration-3D}) when (\ref{Eq:new_time}) is added to
(\ref{Eq:iteration-3D}) and $\phi(t)$ and $\phi(t+1/2)$ are replaced by
$\bar{\phi}(z,v_{3})$ and $\bar{\phi}(z+1/2,v_{3})$, respectively. This
can be
solved as in \cite{GZ15} by taking an ODE with the format $y^{\prime}
=-y^{3}-y$. Hence the system (\ref{Eq:iteration-3D}) must be updated to
\begin{equation}
    \left\{
    \begin{array}
        [c]{l}
        u_{1}^{\prime}=c((\bar{f}_{1}(r(v_{1}),r(v_{2}),r(v_{3}))-u_{1})^{3}+\bar
        {f}_{1}(r(v_{1}),r(v_{2}),r(v_{3}))-u_{1})\bar{\phi}(w,v_{3})\\
        u_{2}^{\prime}=c((\bar{f}_{2}(r(v_{1}),r(v_{2}),r(v_{3}))-u_{2})^{3}+\bar
        {f}_{2}(r(v_{1}),r(v_{2}),r(v_{3}))-u_{2})\bar{\phi}(w,v_{3})\\
        u_{3}^{\prime}=c((\bar{f}_{3}(r(v_{1}),r(v_{2}),r(v_{3}))-u_{3})^{3}+\bar
        {f}_{3}(r(v_{1}),r(v_{2}),r(v_{3}))-u_{3})\bar{\phi}(w,v_{3})\\
        v_{1}^{\prime}=c((r(u_{1})-v_{1})^{3}+r(u_{1})-v_{1})\bar{\phi}(w+1/2,v_{3})\\
        v_{2}^{\prime}=c((r(u_{2})-v_{2})^{3}+r(u_{2})-v_{2})\bar{\phi}(w+1/2,v_{3})\\
        v_{3}^{\prime}=c((r(u_{3})-v_{3})^{3}+r(u_{3})-v_{3})\bar{\phi}(w+1/2,v_{3})\\
        z^{\prime}=1-\zeta_{m-3/16,m-1/8}(v_{3})(z+1).
    \end{array}
    \right.  \label{Eq:iteration-3D-1}
\end{equation}
To show that $x_{halt}$ is a sink of (\ref{Eq:iteration-3D-1}), we
first observe that $x_{halt}$ is an equilibrium point of
(\ref{Eq:iteration-3D-1}). If we are able to show that the Jacobian
matrix $A$ of (\ref{Eq:iteration-3D-1}) at $x_{halt}$ has only
negative eigenvalues, then
$x_{halt}$ will be a sink. A straightforward calculation shows that
\[
A=\left[
\begin{array}
	[c]{cc}
	B & 0\\
	0 & -\zeta_{m-3/16,m-1/8}(m)
\end{array}
\right]
\]
where $B$ is a $6\times 6$ matrix where the entries of the main diagonal take the value $-\bar{\phi}(0,m)$ and all other entries are $0$.
Thus $A$ has two eigenvalues:\ $-\bar{\phi
}(0,m)=-(\phi(0)+1)\leq-1$ and $-\zeta_{m-3/16,m-1/8}(m)=-1$. Since
$A$ only has negative eigenvalues, we conclude that $x_{halt}$ is a
sink of (\ref{Eq:iteration-3D-1}).

We now demonstrate that the basin of attraction of $x_{halt}$
is non-computable. Let $M$ be a universal Turing machine with a
transition function simulated by $\bar{f} = (\bar{f}_{1},
\bar{f}_{2}, \bar{f}_{3})$. Suppose that the initial state of $M$ is
encoded as the number 1 (where the states are encoded as integers
${1,\ldots,m}$ and $m$ is assumed to be the unique halting state).
Then, on input $w$, the initial configuration of $M$ is encoded as
$(0,w,1) \in \mathbb{N}^{3}$. $M$ halts on input $w \in \mathbb{N}$
if and only if $\bar{f}^{[k]}(0,w,1)$ converges to
$\bar{x}_{halt}=(0,0,m)$, and the same is true for any input $x \in
\mathbb{R}^{3}$ satisfying $\left\Vert x-(0,w,1)\right\Vert\leq1/4$.

As shown in the previous section, the basin of attraction of
$\bar{x}_{halt}$ for the discrete dynamical system defined by
$\bar{f}$ cannot be computable. In fact, if the basin of attraction
of $\bar{x}_{halt}$ were computable, then we could solve the Halting
problem as follows: compute a $1/8$-approximation of the basin of
attraction of $\bar{x}_{halt}$. To decide whether $M$ halts with
input $w$, check whether $(0,w,1)$ belongs to that approximation.
Since the halting problem is not computable, the same should be true
for the basin of attraction of $\bar{x}_{halt}$.

We can apply the same idea to ODEs by using the robust iteration of
$\bar{f}$ via the ODE (\ref{Eq:iteration-3D-1}). However, to show a
similar result, we need to prove that any $x\in\mathbb{R}^{7}$
satisfying $\left\Vert x-(0,w,1,0,w,1,0)\right\Vert \leq1/8$ will
converge to $\bar{x}_{halt}$ if and only if $M$ halts with input $w$.
In other words, we need robustness to perturbations in the initial
condition to demonstrate the non-computability of the basin of
attraction of $x_{halt}$, which shows that trajectories starting in a
neighborhood of a configuration encoding an initial configuration
will either all converge to $x_{halt}$ (if $M$ halts with the
corresponding input) or none of these trajectories will converge to
$x_{halt}$ (if $M$ does not halt with the corresponding input).

While the robustness of the convergence to the sink is ensured for
the first six components of $(0,w,1,0,w,1,0)$ due to the robustness
of $\bar{f}$ (at least until $M$ halts), the same does not hold for
the last component $z$, which concerns time. If we start at $t=-1/4$
or $t=1/4$, we begin the periodic cycle required to update the
iteration of $\bar{f}$ too soon or too late. To address this
problem, we modify the function $\phi$ (and thus $\bar{\phi}$ due to
(\ref{Eq:phi_bar})) to ensure that $\phi$ has the additional
property that $\phi(t)=0$ when $t\in\lbrack0,1/4]$, to ensure
robustness to \textquotedblleft late\textquotedblright\ starts (i.e.
when $z\in]0,1/4]$). Note also that $\phi(t)=0$ when
$t\in\lbrack-1/2,0]$, since $z$ is periodic, which ensures
robustness to \textquotedblleft premature\textquotedblright\ starts
(i.e. when $z\in \lbrack-1/4,0[$). Since $\phi$ is periodic with
period 1 and it must be $\phi(t)=0$ when
$t\in\lbrack0,1/4]\cup\lbrack1/2,1]$ and $\phi(t)>0$ when
$t\in]1/4,1/2[$, we take
\[
\phi(t)=\zeta\left(  \sin\left(  2\pi t-\frac{\pi}{4}\right)  -\frac{1}
{\sqrt{2}}\right)  .
\]
Indeed, in the interval $[0,1]$, $\sin\left( 2\pi
t-\frac{\pi}{4}\right) \in\lbrack1/\sqrt{2},1]$ only on $[1/4,1/2]$,
which implies that $\phi(t)=0$ when
$t\in\lbrack0,1/4]\cup\lbrack1/2,1]$ and $\phi(t)>0$ when
$t\in]1/4,1/2[$, due to the properties of $\zeta$. With this
modification, we have ensured robustness to perturbations in the
initial condition for all components of $x$ including time. We can
now conclude, similarly as we did for the map $\bar{f}$, that the
basin of attraction of (\ref{Eq:iteration-3D-1}) must be
non-computable.

This ensures condition (i) and (ii) above.

\subsection{Establishing robust non-computability under perturbations}

In this subsection we improve the construction of the previous subsection to show that condition (iii) also holds, for the conditions presented at the beginning of Section~\ref{Sec:sink}. In other words, we establish that both the sink and the non-computability of the basin of attraction can be made resilient to perturbations

In order to demonstrate that the dynamics of
(\ref{Eq:iteration-3D-1}) remain robust even when subjected to
perturbations, let us consider a function
$g:\mathbb{R}^{7}\rightarrow\mathbb{R}^{7}$ such that $\left\Vert
f-g\right\Vert_{1}\leq1/16$, where (\ref{Eq:iteration-3D-1}) is
expressed as $x^{\prime}=f(x)$. As long as $M$ has not yet halted,
the dynamics of $y^{\prime}=f(x)$ will remain robust against
perturbations to $f$, with the exception of the component $z$ which
is not perturbed. This is because the map
$\bar{f}:\mathbb{R}^{3}\rightarrow\mathbb{R}^{3}$ can robustly
simulate Turing machines, and the dynamics of
(\ref{Eq:iteration-3D-1}) are themselves robust against
perturbations of magnitude $\leq1/16$, as previously demonstrated in
the analysis of (\ref{Eq:robust-iteration}). We should note that we
do not use $\rho=1/8$ as a bound for $\xi(t)$ in
(\ref{Eq:target-perturbed}) since, as previously seen, the total
targeting error $\left\vert x(t)-b\right\vert $ is bounded by
$\gamma+\rho(t_{1}-t_{0})$. However, when $z$ is perturbed, as we
will see, we may not have $t_{1}-t_{0}=1/2$, but instead
$t_{1}-t_{0}\in\lbrack3/4\cdot1/2,5/4\cdot1/2]=[3/8,5/8]$. Using
$\rho=1/16$ instead of $\rho=1/8$ compensates for this issue.

Under these conditions, we can still use $y^{\prime}=g(y)$ to
simulate $M$ until it halts. If we add a perturbation of magnitude
$\leq1/4$ to the right-hand side of the dynamics of $z$ in
(\ref{Eq:iteration-3D-1}), we can conclude that $3/4\leq
z^{\prime}(t)\leq5/4$, meaning that $z(t)$ will remain strictly
increasing and can be used as the ``time variable'' $t$ when iterating
$\bar{f}$. However, there is a potential issue when updating the
iteration cycles of $\bar{f}$ with the ODE
(\ref{Eq:iteration-3D-1}). As previously seen, these cycles occur
over consecutive half-unit time intervals. The issue is that the
first half-unit interval $[0,1/2]$ in a perturbed version of
(\ref{Eq:iteration-3D-1}) will correspond to time values
$t_{1}>t_{0}$ such that $z(t_{0})=0$ and $z(t_{1})=1/2$. Therefore,
when determining the value of $c$ for (\ref{Eq:iteration-3D-1}), we
must use $\int_{t_{0}}^{t_{1}}\phi(z(t))dt$ instead of
$\int_{0}^{1/2}\phi(t)dt$. This will depend on the perturbed value
of $z(t)$, which could potentially lead to issues. However, from
$3/4\leq z^{\prime}(t)\leq5/4$ (which implies that $t_{1}-t_{0}\in
\lbrack3/4\cdot1/2,5/4\cdot1/2]$ as assumed above) we get
$4/3\geq1/z^{\prime
}(t)$, and thus
\[
\int_{t_{0}}^{t_{1}}\phi(z(t))dt=\int_{t_{0}}^{t_{1}}\phi(z(t))z^{\prime
}(t)\frac{1}{z^{\prime}(t)}dt
\]
which implies that, by the change of variables $\tau=z(t)$ (recall that
$\phi(t)\geq0$ for all $t\in\mathbb{R}$)
\[
0<\int_{t_{0}}^{t_{1}}\phi(z(t))dt\leq\frac{4}{3}\int_{0}^{1/2}\phi(\tau)\tau.
\]
Hence, if we take
\begin{equation}
    c\geq\frac{1}{2\gamma^{2}\frac{4}{3}\int_{t_{0}}^{t_{1}}\phi(t)dt}=\frac
    {3}{8\gamma^{2}\int_{t_{0}}^{t_{1}}\phi(t)dt}\label{Eq:c-1}
\end{equation}
we will have enough time to appropriately update each iteration, even if the
\textquotedblleft new\textquotedblright\ time variable $z(t)$ evolves faster
than $t$, thus ensuring robustness to perturbations of the dynamics of
(\ref{Eq:iteration-3D-1}), at least until $M$ halts.

Now let's address the main concern: what happens after $M$ halts. We
will choose $\gamma=1/16$ to ensure that if $M$ halts with input
$w$, then any trajectory of the perturbed system $y^{\prime}=g(y)$
starting in $B(c_{w},1/8)$, where $c_{w}\in\mathbb{N}^{7}$ is the
initial configuration associated with input $w$, will enter
$B(x_{halt},1/4)$ and stay there, where $x_{halt} =(0,0,m,0,0,m,0)$
is the halting configuration. Conversely, if $M$ does not halt with
input $w$, then no trajectory of the perturbed system
$y^{\prime}=g(y)$ starting in $B(c_{w},1/8)$ will enter
$B(x_{halt},1/4)$ (recall that the total error of the perturbed
targeting equation (\ref{Eq:target-perturbed}) is given by
$\left\vert x(t_{1})-b\right\vert <\gamma+\rho(t_{1}-t_{0})$ when
(\ref{Eq:iteration-3D-1}) is actively simulating $M$, i.e., until
$M$ halts).

We first observe that, once the machine $M$ halts at time $t^*$, we
can infer from equation (\ref{Eq:phi_1-halt}) that
$2\geq\bar{\phi}(z(t),v_{3}(t))\geq1$ and
$2\geq\bar{\phi}(-z(t),v_{3}(t))\geq1$. Now, if we rewrite equation
(\ref{Eq:iteration-3D-1}) as $x^{\prime}=f(x)$, we can show that for
any $x\in B(x_{halt},1/4)=\{x:\left\Vert x-x_{halt}\right\Vert
\leq1/4\}$, we have (by using the standard inner product and noticing
the expressions on the right-hand side of (\ref{Eq:iteration-3D-1}))
that:
\[
\left\langle f(x)-x_{halt},x-x_{halt}\right\rangle \leq-c\left\Vert
x-x_{halt}\right\Vert _{2}^{2}\leq-c\left\Vert x-x_{halt}\right\Vert ^{2}.
\]
(Recall also the Euclidean norm $\left\Vert
(x_{1},\ldots,x_{n})\right\Vert
_{2}=\sqrt{x_{1}^{2}+\ldots+x_{n}^{2}}$ for
$(x_{1},\ldots,x_{n})\in\mathbb{R}^{n}$ and that $\left\Vert (x_{1}
,\ldots,x_{n})\right\Vert \leq\left\Vert
(x_{1},\ldots,x_{n})\right\Vert _{2}\leq\sqrt{n}\left\Vert
(x_{1},\ldots,x_{n})\right\Vert $, where $\left\Vert
\cdot\right\Vert $ is the max-norm.) As $c$ must satisfy
(\ref{Eq:c-1}), we can assume without loss of generality
that $c\geq1$, which yields
\begin{equation}
    \left\langle f(x)-x_{halt},x-x_{halt}\right\rangle \leq-\left\Vert
    x-x_{halt}\right\Vert ^{2}\label{Eq:inner-product-h}
\end{equation}
for all $x\in B(x_{halt},1/4)$.

By standard results in dynamical systems (see e.g., \cite[Theorems 1
and 2 of p.~305]{HS74}), there exists some $\varepsilon > 0$ such
that if $\left\Vert g-f\right\Vert_{1}\leq \varepsilon$ (in fact,
this condition only needs to be satisfied on $B(x_{halt},1/4)$), then
$g$ will also have a sink $s_g$ in the interior of
$B(x_{halt},1/16)$. We now assume that $\left\Vert g-f\right\Vert_{1}\leq \min(1/16,\varepsilon)$ on $B(x_{halt},1/4)$.

Next, let us assume that $x \in B(s_g,3/16)$. Since $\left\Vert
s_g-x_{halt}\right\Vert \leq 1/16$, we conclude that $\left\Vert
x-x_{halt}\right\Vert \leq \left\Vert x-s_g\right\Vert +\left\Vert
s_g-x_{halt}\right\Vert \leq 3/16+1/16=1/4$. Therefore, $x\in
B(x_{halt},1/4)$, which implies that (\ref{Eq:inner-product-h})
holds for every $x\in B(s_{g},3/16)$. In what follows, we assume
that $x\in B(s_g,3/16)$. Using (\ref{Eq:inner-product-h}), we
obtain:
\begin{align}
   & \left\langle g(x)-x_{halt},x-s_{g}\right\rangle \nonumber\\
    & =\left\langle
    f(x+x_{halt}-s_{g})-x_{halt},x-s_{g}\right\rangle +\left\langle
    g(x)-f(x+x_{halt}-s_{g}),x-s_{g}\right\rangle \nonumber\\
   &   =\left\langle f(x+x_{halt}-s_{g})-x_{halt},(x+x_{halt}-s_{g})-x_{halt}
    \right\rangle +\left\langle g(x)-f(x+x_{halt}-s_{g}),x-s_{g}\right\rangle
    \nonumber\\
    & \leq-\left\Vert x+x_{halt}-s_{g}-x_{halt}\right\Vert ^{2}+\left\langle
    g(x)-f(x+x_{halt}-s_{g}),x-s_{g}\right\rangle \nonumber\\
    & \leq-\left\Vert x-s_{g}\right\Vert ^{2}+\left\langle g(x)-f(x+x_{halt}
    -s_{g}),x-s_{g}\right\rangle .\label{Eq:inner-product-g}
\end{align}
Furthermore $\alpha(x)=g(x)-f(x+x_{halt}-s_{g})$ is $0$ when $x=s_{g}$ and
\begin{equation}
    \left\Vert D\alpha(x)\right\Vert \leq\left\Vert Dg(x)-Df(x)\right\Vert
    +\left\Vert Df(x)-Df(x+x_{halt}-s_{g})\right\Vert .\label{Eq:Dalpha}
\end{equation}
Since $\left\Vert g-f\right\Vert _{1}\leq\min(1/16,\varepsilon)$ on
$B(x_{halt},1/4)$, this implies that $\left\Vert Dg(x)-Df(x)\right\Vert
\leq1/16$ on $B(x_{halt},1/4)$. Moreover, because $Df$ is continuous on
$B(x_{halt},1/4)$, one can determine some $\delta>0$ such that $\left\Vert
Df(x)-Df(y)\right\Vert \leq1/16$ for all $x,y\in B(x_{halt},1/4)$ satisfying
$\left\Vert x-y\right\Vert \leq\delta$. In particular, if $\left\Vert
x_{halt}-s_{g}\right\Vert \leq\delta$, then (\ref{Eq:Dalpha}) yields
$\left\Vert D\alpha(x)\right\Vert \leq1/16+1/16=1/8$. By classical results
(e.g.$~$\cite[Theorems 1 and 2 of p.~305]{HS74}) we can choose $\varepsilon
_{2}>0$ such that $\left\Vert g-f\right\Vert _{1}\leq\varepsilon_{2}$ implies
$\left\Vert x_{halt}-s_{g}\right\Vert \leq\delta$ as required. Thus when
$\left\Vert g-f\right\Vert _{1}\leq\min\{1/16,\delta,\varepsilon_{2}\}$, we
get that $1/8$ is a Lipschitz constant for $\alpha$ on $B(x_{halt},1/4)$ and
thus
\[
\left\Vert \alpha(x)\right\Vert =\left\Vert \alpha(x)-\alpha(s_{g})\right\Vert
\leq1/8\left\Vert x-s_{g}\right\Vert .
\]
This last inequality and the Cauchy--Schwarz inequality imply that
\begin{align*}
    \left\vert \left\langle g(x)-f(x+x_{halt}-s_{g}),x-s_{g}\right\rangle
    \right\vert  &  =\left\vert \left\langle \alpha(x),x-s_{g}\right\rangle
    \right\vert \\
    &  \leq\left\Vert \alpha(x)\right\Vert _{2}\cdot\left\Vert x-s_{g}\right\Vert
    _{2}\\
    &  \leq7\left\Vert \alpha(x)\right\Vert \cdot\left\Vert x-s_{g}\right\Vert \\
    &  \leq\frac{7}{8}\left\Vert x-s_{g}\right\Vert ^{2}.
\end{align*}
This, together with (\ref{Eq:inner-product-g}), yields
\begin{align*}
    \left\langle g(x),x-s_{g}\right\rangle  &  \leq-\left\Vert x-s_{g}\right\Vert
    ^{2}+\left\langle g(x)-f(x+x_{halt}-s_{g}),x-s_{g}\right\rangle \\
    &  \leq-\left\Vert x-s_{g}\right\Vert ^{2}+\frac{7}{8}\left\Vert
    x-s_{g}\right\Vert ^{2}\\
    &  \leq-\frac{1}{8}\left\Vert x-s_{g}\right\Vert ^{2}
\end{align*}
In particular this shows that $g(x)$ always points inwards inside
$B(s_{g},3/16).$

Since it is well known that
\[
\frac{d}{dt}\left\Vert y(t)\right\Vert _{2}=\frac{1}{\left\Vert
    y(t)\right\Vert _{2}}\left\langle \frac{dy(t)}{dt},y(t)\right\rangle
\]
we get from the last inequality that
\begin{align*}
    \frac{d}{dt}\left\Vert x-s_{g}\right\Vert _{2} &  =\frac{1}{\left\Vert
        x-s_{g}\right\Vert _{2}}\left\langle x^{\prime},x-s_{g}\right\rangle \\
    &  =\frac{1}{\left\Vert x-s_{g}\right\Vert _{2}}\left\langle g(x),x-s_{g}
    \right\rangle \\
    &  \leq\frac{1}{\left\Vert x-s_{g}\right\Vert }\left\langle g(x),x-s_{g}
    \right\rangle \\
    &  \leq-\frac{1}{8}\left\Vert x-s_{g}\right\Vert
\end{align*}
which shows that $\left\Vert x-s_{g}\right\Vert _{2}$ converges exponentially
fast to $s_{g}$ whenever $x\in B\left(  s_{g},3/16\right)  $. Therefore
$B(s_{g},3/16)$ is contained in the basin of attraction of $s_{g}$. In
particular, because $B(x_{halt},1/8)\subseteq B(s_{g},3/16)$, we conclude that
if an initial configuration $c_{w}=(0,w,1,0,w,1,0)\in\mathbb{N}^{7}$ is such
that $M$ halts with input $w$, then a trajectory starting on $B(c_{w},1/4)$ of
the perturbed system $x^{\prime}=g(x)$ of (\ref{Eq:iteration-3D-1}) will reach
$B(x_{halt},1/4)$, and thus $B(s_{g},3/16)$, iff $M$ halts with input $w$.
Furthermore, because any trajectory that enters $B(x_{halt},1/8)\subseteq
B(s_{g},3/16)$ will converge to $s_{g}$ then $M$ halts with input $w$ iff
$B(c_{w},1/4)$ is inside the basin of attraction of $s_{g}$ for $x^{\prime
}=g(x)$ whenever $\left\Vert f-g\right\Vert \leq1/4$ over $\mathbb{R}^{7}$ and
$\left\Vert f-g\right\Vert _{1}\leq\min\{1/16,\delta,\varepsilon_{2}\}$ over
$B(x_{halt},1/4)$. Indeed, if $M$ does not halt with input $w$, then any
trajectory which starts on $B(c_{w},1/4)$ will never enter $B(x_{halt},1/4)$
under the dynamics of $x^{\prime}=g(x)$ and thus never enter $B(s_{g},3/16)$,
otherwise it would converge to $s_{g}$ and then enter $B(s_{g},1/16)\subseteq
B(x_{halt},1/4)$, a contradiction. Using similar arguments to those used for
$f$, we conclude that the basin of attraction for $g$ is not computable. This ends the proof of Theorem A.\medskip

We briefly mention that in the context of the continuous dynamical
system $y' = f(y)$, the function $f$ is $C^\infty$ (infinitely
differentiable) rather than analytic, as is the case in the discrete
counterpart. The absence of analyticity in $f$ stems from the
function $\phi$ employed to construct it (recall that $\phi(t)=0$ on
the intervals $(k, k+\frac{1}{2})$ for integers $k$). However, by
employing a more sophisticated $\phi$ as described in \cite{GZ15},
it becomes possible to enhance $f$ to an analytic function. For the
sake of readability, we have chosen to present an example of a
$C^{\infty}$ system.

\section{Proof of Theorem C -- Basins of attraction of structurally\texorpdfstring{\\}{}
stable planar systems are uniformly computable}

\label{Sec:RobustNoncomputability}

In the previous section, we demonstrated the existence of a
$C^{\infty}$ and computable system (\ref{eq_main}) that possesses a
computable sink with a non-computable basin of attraction. Moreover,
this non-computability persists throughout a neighborhood of $f$. It
should be noted that a dynamical system is locally stable near a
sink. Thus our example shows that local stability at a sink does not
guarantee the existence of a numerical algorithm capable of
computing its basin of attraction.

In this section, we investigate the relationship between the global
stability of a planar structurally stable system (\ref{eq_main}) defined over the unit ball and
the computability of its basins of attraction. We demonstrate that
if the system is globally stable, then the basins of attraction of
all of its sinks are computable. This result highlights that global
stability is not only a strong analytical property but also gives
rise to strong computability regarding the computation of basins of
attraction. Moreover, it shows that strong computability is
``typical'' on compact planar systems since it is well known (see
e.g.~\cite[Theorem 3 on p.~325]{Per01}) that in this case the set of
$C^1$ structurally stable vector fields is open and dense over the
set of $C^1$ vector fields.

We begin this section by introducing some preliminary definitions.
Let $K$ be a closed disk in $\mathbb{R}^2$ centered at the origin
with a rational radius. In particular, let $\mathbb{D}$ denote the
closed unit disk of $\mathbb{R}^2$. We define $\mathcal{V}(K)$ to be
the set of all $C^1$ vector fields mapping $K$ to $\mathbb{R}^2$
that point inwards along the boundary of $K$. Furthermore, we define
$\mathcal{O}_2$ to be the set of all open subsets of $\mathbb{R}^2$
equipped with the topology generated by the open rational disks,
i.e., disks with rational centers and rational radii, as a subbase.

For a structurally stable planar system $x^{\prime}=f(x)$ defined on
the closed disk $K$, it has only finitely many equilibrium points
and periodic orbits, and all of them are hyperbolic (see \cite{Pei59}). Recall from
Section~\ref{Subsec:DynamicalSystem} that a point $x_{0} \in K$ is
called an equilibrium point of the system if $f(x)=0$, since any
trajectory starting at an equilibrium stays there for all
$t\in\mathbb{R}$. Recall also that an equilibrium point $x_0$ is called hyperbolic if
all the eigenvalues of $Df(x_{0})$ have non-zero real parts. If both
eigenvalues of $Df(x_{0})$ have negative real parts, then it can be
shown that $x_{0}$ is a sink. A sink attracts nearby trajectories.
If both eigenvalues have positive real parts, then $x_{0}$ is called
a source. A source repels nearby trajectories. If the real parts of
the eigenvalues have opposite signs, then $x_{0}$ is called a
saddle (see Figure \ref{Fig:basin} for a picture of a saddle point). A saddle attracts some points (those lying in the stable
manifold, which is a one-dimensional manifold for the planar
systems), repels other points (those lying in the unstable manifold,
which is also a one-dimensional manifold for the planar systems,
transversal to the stable manifold), and all trajectories starting
in a neighborhood of a saddle point but not lying on the stable
manifold will eventually leave this neighborhood. A periodic orbit
(or limit cycle) is a closed curve $\gamma$ with the property that
there is some $T>0$ such that $\phi(f, x)(T) = x$ for any
$x\in\gamma$. Hyperbolic periodic orbits have properties similar to
hyperbolic equilibria. For a planar system, there are only
attracting or repelling hyperbolic periodic orbits. See \cite[p.
225]{Per01} for more details.

In this section, we demonstrate the existence of an algorithm that
can compute the basins of attraction of sinks for any structurally
stable planar vector field defined on a compact disk $K$ of
$\mathbb{R}^2$. Furthermore, this computation is uniform across the
entire set of such vector fields.

In Theorem C
below, we consider the case
where $K=\mathbb{D}$ for simplicity, but the same argument applies
to any closed disk with a rational radius. Before stating and
proving Theorem C, we present two lemmas,
the proofs of which can be found in \cite{GZ21}. Let $\mathcal{SS}_2
\subset \mathcal{V}(\mathbb{D})$ be the set of all $C^1$
structurally stable planar vector fields defined on $\mathbb{D}$.

\begin{lem} \label{Lem:NumberSinks} The map $\Psi_{N}: \mathcal{SS}_2
\to \mathbb{N}$, $f\mapsto \Psi_{N}(f)$, is computable, where
$\Psi_{N}(f)$ is the number of the sinks of $f$ in $\mathbb{D}$.
\end{lem}

\begin{lem} The map $\Psi_{S}: \mathcal{SS}_2 \times \mathbb{N} \to \mathcal{P}(\mathbb{R}^2)$ is computable, where $\Psi_{S}(f,k)$ returns a set of disjoint $1/n\times1/n$ squares, where $n\in\mathbb{N}$ is such that $n\geq k$ and each square is centered at a rational point. Furthermore, each square has exactly one equilibrium point (zero) of $f$.
\end{lem}

\begin{mainthm}{3}
The map $\Psi: \mathcal{SS}_2 \times \mathbb{D} \to \mathcal{O}$ is computable, where
$\Psi(f,s)=W_s$ is the basin of attraction of the sink $s$.
\end{mainthm}

\begin{proof} Let us fix an $f\in \mathcal{SS}_2$. Assume that $\Psi_{N}(f)\neq 0$
and $s$ is a sink of $f$. In \cite{Zho09} and \cite{GZ22}, it has
been shown that:
\begin{enumerate}
\item $W_{s}$ is a r.e.~open subset of $\mathbb{D}\subseteq\mathbb{R}^2$;
\item there is an algorithm that on input $f$ and $k\in \mathbb{N}$, $k>0$, computes a finite
sequence of mutually disjoint closed squares or closed ring-shaped
strips (annulus) such that (see Figure \ref{fig:rings}):
\begin{enumerate}
\item each square contains
exactly one equilibrium point with a marker indicating if it
contains a sink, a source, or a saddle;
\item each annulus contains exactly one
periodic orbit with a marker indicating if it contains an attracting
or a repelling periodic orbit;
\item each square (resp. annulus) containing a sink (resp. an attracting periodic orbit) is
time invariant for $t\geq0$;
\item the union of this finite sequence contains all equilibrium
points and periodic orbits of $f$, and the Hausdorff distance
between this union and the set of all equilibrium points and
periodic orbits is less than $1/k$;
\item for each annulus, $1\leq i\leq p(f)$, the minimal distance
between the inner boundary (denoted as $IB_i$) and the outer
boundary (denoted as $OB_i$), $m_i = \min\{ d(x, y) : \, x\in IB_i,
y\in OB_i\}$, is computable from $f$ and $m_i>0$.
\end{enumerate}
\end{enumerate}

\begin{figure}
	\begin{center}
		\includegraphics[width=7cm]{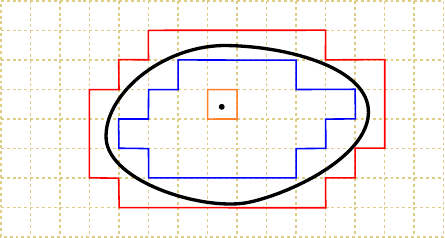}
	\end{center}
	\caption{Result of the algorithm from \cite{GZ22} which computes hyperbolic equilibrium points and hyperbolic periodic orbits with some given (input) accuracy. The periodic orbit is surrounded by a (red) outer boundary and an inner (blue) boundary which delimitates a region approximating the periodic orbit. The orange square delimitates an equilibrium point.}\label{fig:rings}
\end{figure}

We begin with the case that $f$ has no saddle point. Since $W_s$ is
r.e.~open, there exists computable sequences $\{ a_n\}$ and $\{
r_n\}$, $a_n\in \mathbb{Q}^2$ and $r_n\in \mathbb{Q}$, such that
$W_s=\cup_{n=1}^{\infty}B(a_n, r_n)$.

Let $A$ be the union of all
squares and annuli in the finite sequence containing a sink or an
attracting periodic orbit except the square containing $s$, and let
$B$ be the union of all sources and repelling periodic orbits. Note that a source is an equilibrium point (even if unstable) and thus will not belong to $W_s$. Similarly each repelling periodic orbit is an invariant set and thus will also not belong to $W_s$. Periodic orbits and equilibrium points are closed sets and thus $B$ is a closed set of $\mathbb{D}$, which is also computable due to the results from \cite{GZ22} mentioned above.
Hence, $\mathbb{D}\setminus B$ is a computable open subset of
$\mathbb{D}$. Moreover, since $f$ has no saddle, $W_s\subset
\mathbb{D}\setminus B$. List the squares in $A$ as $S_1, \ldots,
S_{e(f)}$ and annuli as $C_1, \ldots, C_{p(f)}$. Denote the center
and the side-length of $S_j$ as $CS_j$ and $l_j$, respectively, for
each $1\leq j\leq e(f)$.

We first present an algorithm -- the classification algorithm --
that for each $x\in\mathbb{D}\setminus B$ determines whether $x\in W_s$ or $x$ is in the union of basins
of attraction of the sinks and attracting periodic orbits contained
in $A$. The algorithm works as follows: for each $x\in
\mathbb{D}\setminus B$, simultaneously compute

\[ \left\{ \begin{array}{l}
   d(x, a_n), n=1, 2, \ldots \\ \\
   d(\phi_{t}(x), CS_j), 1\leq j\leq e(f), t=1, 2, \ldots \\ \\
   \mbox{$d(\phi_{t}(x), IB_i)$ and $d(\phi_{t}(x), OB_i)$}, 1\leq
   i \leq p(f), t=1, 2, \ldots \end{array} \right. \]
where $\phi_t(x)=\phi(f, x)(t)$ is the solution of the system $dz/dt
= f(z)$ with the initial condition $z(0)=x$ at time $t$. (Recall
that the solution, as a function of time $t$, of the initial-value
problem is uniformly computable from $f$ and $x$ \cite{GZB07}.) Halt
the computation whenever one of the following occurs: (i) $d(x,
a_n)<r_n$; (ii) $d(\phi_{t}(x), CS_j) < l_j/2$ for some $t=l\in
\mathbb{N}$ ($l>0$); or (iii) $d(\phi_{t}(x), IB_i) < m_i$ and
$d(\phi_{t}(x), OB_i) < m_i$ for $t=l\in \mathbb{N}$ ($l>0$). If the
computation halts, then either $x\in W_{s}$ provided that $d(x,
a_n)<r_n$ or else $\phi_{t}(x)\in S_j$ or $\phi_{t}(x)\in C_i$ for
some $t=l>0$. Since $S_j$ and $C_i$ are time invariant for $t>0$ (this follows from the results of \cite{GZ22}),
each $S_j$ contains exactly one sink for $1\leq j\leq e(f)$, and
each $C_i$ contains exactly one attracting periodic orbit for $1\leq
i\leq p(f)$, it follows that either $x$ is in the basin of
attraction of the sink contained in $S_j$ if (ii) occurs or $x$ is
in the basin of attraction of the attracting periodic orbit
contained in $C_i$ if (iii) occurs. We note that, for any $x\in
\mathbb{D}\setminus B$, exactly one of the halting status, (i),
(ii), or (iii), can occur following the definition of $W_s$ and the
fact that $S_j$ and $C_i$ are time invariant for $t>0$. Let $W_{A}$
be the set of all $x\in \mathbb{D}\setminus B$ such that the
computation halts with halting status (ii) or (iii) on input $x$.
Then it is clear that $W_{s}\cap W_{A}=\emptyset$.

We turn now to show that the computation will halt. Since there is
no saddle, every point of $\mathbb{D}$ that is not a source or on a
repelling periodic orbit will either be in $W_s$ or the trajectory starting on that point will converge to a
sink/attracting periodic orbit contained in $A$ as $t\to \infty$
(this is ensured by the structural stability of the system and
Peixoto's characterization theorem; see, for example, \cite{Pei59}).

Thus either $x\in W_s$ or $x$ will eventually enter some $S_j$ (or
$C_i$) and stay there afterwards for some sufficiently large
positive time $t$. Hence the condition (i) or (ii) or (iii) will be
met for some $t>0$.

Since $W_s$ is a r.e.~open set due to the results of \cite{Zho09}, to prove that $W_s$ is computable it is suffices to show that the
closed subset $\mathbb{D}\setminus W_{s} = W_{A}\cup B$ is
r.e.~closed; or, equivalently, $W_{A}\cup B$ contains a computable
sequence that is dense in $W_{A}\cup B$ (see e.g.~\cite[Proposition 5.12]{BHW08}). To see this, we first note
that $\mathbb{D}\setminus B$ has a computable sequence as a dense
subset. Indeed, since $\mathbb{D}\setminus B$ is computable open,
there exist computable sequences $\{ z_i\}$ and $\{\theta_i\}$,
$z_i\in \mathbb{Q}^2$ and $\theta_i\in \mathbb{Q}$, such that
$\mathbb{D}\setminus B = \cup_{i=1}^{\infty}B(z_i, \theta_i)$. Let
$\mathcal{G}_{l}=\{ (m/2^{l}, n/2^{l}): \, \mbox{$m, n$ are integers
and $-2^{l}\leq m, n\leq 2^{l}$}\}$ be the $\frac{1}{2^l}$-grid on
$\mathbb{D}$, $l\in \mathbb{N}$. The following procedure produces a
computable dense sequence of $\mathbb{D}\setminus B$: For each input
$l\in \mathbb{N}$, compute $d(x, z_i)$, where $x\in \mathcal{G}_l$
and $1\leq i\leq l$ and output those $\frac{1}{2^l}$-grid points $x$
if $d(x, z_i)<\theta_i$ for some $1\leq i\leq l$. By a standard
paring, the outputs of the computation form a computable dense
sequence, $\{ q_i\}_{i\in\mathbb{N}}$,  of $\mathbb{D}\setminus B$. We now want to obtain a computable dense sequence in $W_A$. If we are able to show that such a computable sequence exists, then it follows that $W_{A}\cup B$ contains
a computable dense sequence. The conclusion comes from the fact that
$B$ is a computable closed subset; hence $B$ contains a computable
dense sequence.

 Then
using the previous classification algorithm one can enlist those points in the sequence $\{ q_i\}_{i\in\mathbb{N}}$ which fall inside
$W_{A}$, say $\tilde{q}_1, \tilde{q}_2, \ldots$. Clearly, $\{
\tilde{q}_j\}_{j\in\mathbb{N}}$ is a computable sequence.

It remains to show that $\{ \tilde{q}_j\}$ is dense in $W_{A}$. It
suffices to show that, for any $x\in W_{A}$ and any neighborhood
$B(x, \epsilon)\cap W_{A}$ of $x$ in $W_{A}$, there exists some
$\tilde{q}_{j_0}$ such that $\tilde{q}_{j_0}\in B(x, \epsilon)\cap
W_{A}$, where $\epsilon>0$ and the disk $B(x, \epsilon)\subset
\mathbb{D}\setminus B$. We begin by recalling a well-known fact that
the solution $\phi_{t}(x)$ of the initial value problem $dx/dt =
f(x)$, $\phi_{0}(x)=x$, is continuous in time $t$ and in initial
condition $x$. In particular, the following estimate holds true for
any time $t>0$ (see e.g.~\cite{BR89}):
\begin{equation} \label{error-initial-condition}
\|\phi_{t}(x) - \phi_{t}(y)\| \leq \| x - y\|e^{Lt}
\end{equation}
where $x=\phi_{0}(x)$ and $y=\phi_{0}(y)$ are initial conditions,
and $L$ is a Lipschitz constant satisfied by $f$. (Since $f$ is
$C^1$ on $\mathbb{D}$, it satisfies a Lipschitz condition and a
Lipschitz constant can be computed from $f$ and $Df$.) Since $x\in
W_{A}$, the halting status on $x$ is either (ii) or (iii). Without
loss of generality we assume that the halting status of $x$ is (ii).
A similar argument works for the case where the halting status of
$x$ is (iii).  It follows from the assumption that $d(\phi_{t}(x),
S_j)<l_j/2$ for some $1\leq j\leq e(f)$ and some $t=l>0$. Compute a
rational number $\alpha$ satisfying $0<\alpha < l_j/2 -
d(\phi_{t}(x), S_j)$ and compute another rational number $\beta$
such that $0< \beta < \epsilon$ and $\| y_1 - y_2\|e^{l\cdot L}<
\alpha$ whenever $\| y_1 - y_2\|<\beta$. Then for any $y\in B(x,
\beta)$,
\begin{eqnarray*}
& & d(\phi_{t}(y), S_j) \\
& \leq & d(\phi_{t}(y), \phi_{t}(x)) + d(\phi_{t}(x), S_j)
\\
& \leq & \alpha + d(\phi_{t}(x), S_j) < (l_j/2) - d(\phi_{t}(x),
S_j) + d(\phi_{t}(x), S_j) = l_j/2
\end{eqnarray*}
which implies that $B(x, \beta)\subset W_{A}$. Since $B(x,
\beta)\subset B(x, \epsilon) \subset \mathbb{D}\setminus B$ and $\{
q_i\}$ is dense in $\mathbb{D}\setminus B$, there exists some
$q_{i_0}$ such that $q_{i_0}\in B(x, \beta)$. Since $B(x,
\beta)\subset W_{A}$, it follows that $q_{i_0}=\tilde{q}_{j_0}$ for
some $j_0$. This shows that $\tilde{q}_{j_0}\in B(x, \epsilon)\cap
W_A$.

We turn now to the general case where saddle point(s) is present. We
continue using the notations introduced for the special case where
the system has no saddle point. Assume that the system has the
saddle points $d_{m}$, $1\leq m\leq d(f)$ and $D_m$ is a closed
square containing $d_m$, $1\leq m\leq d(f)$. For any given $k\in
\mathbb{N}$ ($k>0$), the algorithm constructed in \cite{GZ21} will
output $S_j$, $C_i$, and $D_m$ such that each contains exactly one
equilibrium point or exactly one periodic orbit, the (rational)
closed squares and (rational) closed annuli are mutually disjoint,
each square has side-length less than $1/k$, and the Hausdorff
distance between $C_i$ and the periodic orbit contained inside $C_i$
is less than $1/k$, where $1\leq j\leq e(f)$, $1\leq m\leq d(f)$,
and $1\leq i\leq p(f)$. For each saddle point $d_m$, it is proved in
\cite{GZB12} that the stable manifold of $d_m$ is locally computable
from $f$ and $d_m$; that is, there is a Turing algorithm that
computes a bounded curve -- the flow is planar and so the stable
manifold is one dimensional -- passing through $d_m$ such that
$\lim_{t\to \infty}\phi_{t}(x_0)=d_m$ for every $x_0$ on the curve.
In particular, the algorithm produces a computable dense sequence on
the curve. Pick two points, $z_1$ and $z_2$, on the curve such that
$d_m$ lies on the segment of the curve from $z_1$ to $z_2$. Since
the system is structurally stable, there is no saddle connection;
i.e.~the stable manifold of a saddle point cannot intersect the
unstable manifold of the same saddle point or of another saddle
point. Thus, $\phi_{t}(z_1)$ and $\phi_{t}(z_2)$ will enter $C_B$
for all $t\leq -T$ for some $T>0$, where $C_{B}=(\cup \{ \overline{S}_j : \,
s_j\in B\}) \cup (\cup \{ \overline{C}_i: \, p_i\subset B\})$, where $\overline{S}_j$ and $\overline{C}_i$ denote the squares and annuli computed by the algorithm of \cite{GZ22} which contain repelling equilibrium points (sources) and repelling periodic orbits, respectively. We denote the
curve $\{ \phi_{t}(z_1): \, -T \leq t\leq 0\} \cup \{ z: \mbox{$z$
is on the stable manifold of $d_m$ between $z_1$ and $z_2$}\} \cup
\{ \phi_{t}(z_2): \, -T \leq t\leq 0\}$ as $\Gamma_{d_m}$. Let
$\widetilde{C}=C_B\cup \{ \Gamma_{d_m}: \, 1\leq m\leq d(f)\}$. Then
$\widetilde{C}$ is a computable compact subset in $\mathbb{D}$.
Moreover, every point in $\mathbb{D}\setminus \widetilde{C}$
converges to either a sink or an attracting periodic orbit because
there is no saddle connection. Using the classification algorithm
and a similar argument as above we can show that $W_{A}\cap
(\mathbb{D}\setminus \widetilde{C})$ is a computable open subset in
$\mathbb{D}\setminus \widetilde{C}$ and thus computable open in
$\mathbb{D}$ because $W_{A}\subset
(\mathbb{D}\setminus \widetilde{C})$. Since $W_{A}\subset \mathbb{D}\setminus B$ and
$W_{A}\cap \Gamma_{d_m}=\emptyset$, it follows that
\begin{eqnarray*}
& & d_{H}\left(\mathbb{D}\setminus (W_{A}\cap (\mathbb{D}\setminus
\widetilde{C})), \, \mathbb{D}\setminus (W_{A}\cap
(\mathbb{D}\setminus B))\right) \\
& = & d_{H}\left((\mathbb{D}\setminus W_A)\cup C_{B}, \,
(\mathbb{D}\setminus W_A)\cup B\right) \\
& \leq & d_{H}(C_{B}, B) < \frac{1}{k}.
\end{eqnarray*}
We have thus proved that there is an algorithm that, for each input $k\in
\mathbb{N}$ ($k>0$), computes an open subset $U_{k}=W_{A}\cap
(\mathbb{D}\setminus \widetilde{C})$ of $\mathbb{D}$ such that
$U_{k}\subset W_{A}$ and $d_{H}(\mathbb{D}\setminus U_{k}, \,
\mathbb{D}\setminus W_{A})<\frac{1}{k}$. This shows that $W_{A}$ is
a computable open subset of $\mathbb{D}$. (Recall an equivalent
definition for a computable open subset of $\mathbb{D}$: an open
subset $U$ of $\mathbb{D}$ is computable if there exists a sequence
of computable open subsets $U_k$ of $\mathbb{D}$ such that $U=\cup
U_{k}$ and $d_{H}(\mathbb{D}\setminus U_k, \, \mathbb{D}\setminus
U)\leq \frac{1}{k}$ for every $k\in \mathbb{N}\setminus \{ 0\}$.)
\end{proof}

\begin{cor} For every $f\in \mathcal{SS}_2$ there is a
neighborhood of $f$ in $C^1(\mathbb{D})$ such that the function
$\Psi$ is (uniformly) computable in this neighborhood.
\end{cor}

\begin{proof} The corollary follows from Peixoto's density theorem
and Theorem C.
\end{proof}

\bibliography{ContComp}
\bibliographystyle{alphaurl}

\end{document}